\documentclass{article}
\input epsf.tex
\usepackage{latexsym}
\usepackage{amsfonts,amssymb}
\usepackage{epsfig}
\usepackage{psfig}   

\newtheorem{theorem}{Theorem}[section]
\newtheorem{lemma}[theorem]{Lemma}
\newtheorem{corollary}[theorem]{Corollary}
\newtheorem{remark}[theorem]{Remark}
\newtheorem{remarks}[theorem]{Remarks}
\newtheorem{question}[theorem]{Question}
\newtheorem{problem}[theorem]{Problem}
\newtheorem{definition}[theorem]{Definition}

\newtheorem{proposition}[theorem]{Proposition}

\newtheorem{example}[theorem]{Example}

\title{Equivalent curves in surfaces}
\author{Christopher J. Leininger}
\begin{document}

\maketitle

\section{Introduction}

Given a closed orientable surface $S$ of genus $g > 1$, there exists distinct homotopy classes of curves $\gamma$ and $\gamma'$ on $S$ such that
$$length_{X}(\gamma) = length_{X}(\gamma')$$
for every hyperbolic structure $X$ on $S$ (see \cite{H} and \cite{Ran}).
The proof involves writing down words in the fundamental group which are not conjugate, but because of certain trace relations, they have the same trace squared with respect to every representation into $PSL_{2}({\mathbb R})$ (see Section \ref{traces}).

The hyperbolic metrics are of course very special, and it seems reasonable that if $\gamma$ and $\gamma'$ are as above, there should be some purely topological reason for this.

\begin{problem} \label{basmaj} (\cite{Bas})
Find a topological description of pairs of distinct homotopy classes of curves $\gamma$ and $\gamma'$ which have the same length with respect to every hyperbolic structure on $S$.
\end{problem}

More generally, given a family of path metrics ${\cal G}(S)$ on $S$, we ask the following.

\begin{problem} \label{lenquest}
Do there exist pairs of distinct homotopy classes of curves $\gamma$ and $\gamma'$ which have the same length with respect to every metric in ${\cal G}(S)$?
If so, find a topological description of such pairs.
\end{problem}

For a generic Riemannian metric on $S$, there are {\em no} homotopy classes of curves which have the same length (see \cite{Abr}), so that the answer to the question will {\em generically} be ``no''.
Even when the answer is ``yes'', a solution to the problem is likely to be uninteresting if ${\cal G}(S)$ is an arbitrary family of metrics.
We are thus primarily interested in this problem in the situation that ${\cal G}(S)$ arises {\em naturally}.

The hyperbolic metrics arise (for example) as the unique representatives of their conformal classes with constant curvature $-1$, and are often useful in considering topological problems about surfaces (see e.g. \cite{Tdiff}).
This is probably the most interesting setting for Problem \ref{lenquest}.
Another family of metrics which occur naturally (especially in Teichm\"uller theory) are those metric structures induced by quadratic differentials on $S$, holomorphic with respect to a complex structure (see e.g. \cite{Ab}).
These are Euclidean cone metrics with some additional structure (see Section \ref{quaddiff}).
We will refer to such metrics as {\em branched flat metrics}.

We discuss Problem \ref{lenquest} in the context of both hyperbolic and branched flat metrics, relate the two, and give a complete solution in the latter case.
We also show that an ``obvious'' necessary topological condition with respect to the hyperbolic metrics is not sufficient, and mention a possible alternative.

To state the main theorem, we make the following definitions (see Sections \ref{background} and \ref{quaddiff} for more details).
\begin{definition}
Let $\gamma$ and $\gamma'$ be homotopy classes of essential closed curves on $S$.
\begin{itemize}
\item $\gamma$ and $\gamma'$ are {\em hyperbolically equivalent}, $\gamma \equiv_{h} \gamma'$, if for every hyperbolic structure $X$ on $S$, $length_{X}(\gamma) = length_{X}(\gamma')$.
\item $\gamma$ and $\gamma'$ are {\em branched flat equivalent}, $\gamma \equiv_{bf} \gamma'$, if for every branch flat structure $Q$ on $S$, $length_{Q}(\gamma) = length_{Q}(\gamma')$.
\item $\gamma$ and $\gamma'$ are {\em simple intersection equivalent}, $\gamma \equiv_{si} \gamma'$, if for every isotopy class of essential {\em simple} closed curve $\alpha$ on $S$, $i(\alpha,\gamma) = i(\alpha,\gamma')$.
\item $\gamma$ and $\gamma'$ are {\em trace equivalent}, $\gamma \equiv_{tr} \gamma'$, if for every representation $\rho : \pi_{1}(S) \rightarrow SL_{2}({\mathbb C})$, $tr(\rho(\hat{\gamma}))^{2} = tr(\rho(\hat{\gamma}'))^{2}$, where $\hat{\gamma}$ and $\hat{\gamma}'$ are elements of $\pi_{1}(S)$ representing the conjugacy classes determined by $\gamma$ and $\gamma'$, respectively.
\end{itemize}
The respective equivalence classes will be called {\em hyperbolic classes}, {\em branched flat classes}, {\em simple intersection classes}, and {\em trace classes}.
\end{definition}

The main theorem is the following.

\begin{theorem} \label{main}
Given a pair of homotopy classes of essential closed curves $\gamma$ and $\gamma'$ on $S$, then
$$ \gamma \equiv_{tr} \gamma' \Leftrightarrow \gamma \equiv_{h} \gamma' \Rightarrow \gamma \equiv_{si} \gamma' \Leftrightarrow \gamma \equiv_{bf} \gamma'$$

Moreover, we can find $\gamma$ and $\gamma'$, each of which fill the surface, so that $\gamma \equiv_{si} \gamma'$ yet $\gamma \not\equiv_{h} \gamma'$.
\end{theorem}

We say that a homotopy class of essential closed curves {\em fills} $S$ if every representative intersects every other essential closed curve.

The outline of the paper is as follows.
In Section 2, we give a few standard definitions to fix notation.
In Section 3, we briefly discuss the construction of nontrivial pairs $\gamma$ and $\gamma'$ with $\gamma \equiv_{h} \gamma'$, following Horowitz and Randol.
We then prove that $\gamma \equiv_{tr} \gamma' \Leftrightarrow \gamma \equiv_{h} \gamma'$.
One direction is immediate from the trace formula for lengths, and the other direction follows easily from the main result of \cite{Rap}.
In Section 4, we prove that $\gamma \equiv_{si} \gamma' \Leftrightarrow \gamma \equiv_{bf} \gamma'$.
In the process, we determine an interesting relationship between lengths and heights (see Lemma \ref{lengthfromheight}).
In section 5, we discuss Bonahon's interpretation of Thurston's compactification of Teichm\"uller space which easily implies $\gamma \equiv_{h} \gamma' \Rightarrow \gamma \equiv_{si} \gamma'$.
In Section 6, we construct the required counterexamples, completing the proof of Theorem \ref{main}.
We also describe a few of the difficulties encountered in finding such pairs of curves.
In Section 7, we describe a candidate for a topological characterizations of hyperbolic equivalence and state a couple questions.
We end by citing a few references for more information on Problem \ref{basmaj} and related problems.\\

\noindent
{\bf Acknowledgment:} I would like to thank Jim Anderson, Andrew Casson and my advisor Alan Reid for useful conversations.

\section{Background}  \label{background}

In this section we recall a few standard definitions and fix some notation.

\subsection{Teichm\"uller space and closed curves}  \label{teichprelim}

Let $S$ denote a closed orientable surface of genus $g >1$.
Let ${\cal T}eich(S)$ denote the Teichm\"uller space of $S$, which we think of as the space of hyperbolic metrics on $S$ up to isotopy.
By the Classical Uniformization Theorem (see \cite{FK}), this is equivalently the space of complex structures on $S$ up to isotopy.

We let ${\cal C}rv(S)$ denote the set of homotopy classes of (unoriented) primitive essential closed curves on $S$.
We denote the set of elements in ${\cal C}rv(S)$ with {\em simple} representatives by ${\cal S}(S)$.
Of course, ${\cal S}(S)$ can be naturally identified with the set of isotopy classes of essential simple closed curves.

The geometric intersection number of two elements $\gamma, \gamma' \in {\cal C}rv(S)$ is defined by
$$i(\gamma,\gamma')= \min_{\overline{\gamma} \in \gamma \, , \, \overline{\gamma}' \in \gamma'} |(\overline{\gamma} \times \overline{\gamma}')^{-1}(\Delta)|.$$
where $\overline{\gamma} \times \overline{\gamma}' : S^{1} \times S^{1} \rightarrow S \times S$, $\Delta \subset S \times S$ is the diagonal, and $|A|$ denotes the cardinality of a set $A$.

Given an element $\gamma \in {\cal C}rv(S)$, we also have its associated length function
$$length_{(\cdot)}(\gamma) : {\cal T}eich(S) \rightarrow {\mathbb R}$$
defined, for every hyperbolic metric $X \in {\cal T}eich(S)$, by
$$length_{X}(\gamma) = \inf_{\overline{\gamma} \in \gamma} length_{X}(\overline{\gamma})$$
This is also the length of the geodesic representative of $\gamma$ in $X$.

We will often find it useful to fix a hyperbolic metric on $S$ for reference.
Therefore, whenever we refer to geodesics on $S$ (or in its universal cover) this is meant with respect to the fixed hyperbolic structure, unless otherwise stated.
Once this is done, each element of ${\cal C}rv(S)$ has a unique geodesic representative, and we can identify ${\cal C}rv(S)$ with the set of closed (primitive) geodesics on $S$.
Moreover, geometric intersection number is realized by these representatives (see e.g. \cite{CB} for a proof in the case of simple curves).

\subsection{Measured laminations and foliations} \label{lamprelim}

We denote the space of measured geodesic laminations on $S$ by ${\cal ML}(S)$.
A point $(\Lambda,\lambda) \in {\cal ML}(S)$ consists of a geodesic lamination $\Lambda$ along with a transverse measure $\lambda$ of full support (when no confusion arises, we refer to $(\Lambda,\lambda)$ simply as $\lambda$).
If $\alpha \in {\cal S}(S)$ and $t \in {\mathbb R}_{+}$, then the geodesic representative for $\alpha$, along with $t$ times the transverse counting measure defines a point in ${\cal ML}(S)$.
This defines an injection ${\mathbb R}_{+} \times {\cal S}(S)$ into ${\cal ML}(S)$, and we identify ${\mathbb R}_{+} \times {\cal S}(S)$ with its image.
The relevant theorem here is the following (\cite{Tnotes}, \cite{Bbouts}; see also Section \ref{geodesiccurrents}).

\begin{theorem} \label{mlextend}
(Thurston) ${\mathbb R}_{+} \times {\cal S}(S)$ is dense in ${\cal ML}(S)$.
If $\gamma \in {\cal C}rv(S)$, then
$$(t,\alpha) \rightarrow t \cdot i(\alpha,\gamma)$$
extends uniquely to a continuous function
$$i(\cdot,\gamma) : {\cal ML}(S) \rightarrow {\mathbb R}$$
\end{theorem}
\rightline{$\Box$}

\begin{remarks} \

\noindent
1. This is not exactly the statement in \cite{Tnotes}, but the work there provides this version easily.
This is also a direct consequence of the work in \cite{Bbouts} and \cite{BTeich}.\\
2. The extension of $i(\cdot, \gamma)$ can be explicitly described by
$$ i(\lambda,\gamma) = \inf_{\overline{\gamma} \in \gamma} \int_{\overline{\gamma}} d \lambda $$
\end{remarks}

The space ${\cal ML}(S)$ is naturally homeomorphic to ${\cal MF}(S)$, the space of measure classes of measured singular foliations on $S$.
A point $({\cal F},\mu) \in {\cal MF}(S)$ consists of a singular foliation ${\cal F}$ and a transverse measure $\mu$ of full support (again, we write $\mu=({\cal F},\mu)$ when no confusion arises).

The homeomorphism
$$\Phi:{\cal ML}(S) \rightarrow {\cal MF}(S)$$
is natural in the following sense \cite{L}.

\begin{theorem} \label{mlmf}
(Thurston) For any $\gamma \in {\cal C}rv(S)$, we have
$$i(\lambda,\gamma)
= \inf_{\overline{\gamma} \in \gamma} \int_{\overline{\gamma}} d \Phi(\lambda).$$
\end{theorem}
\rightline{$\Box$}

For obvious reasons, given $\gamma \in {\cal C}rv(S)$ and $\mu \in {\cal MF}(S)$, we write
$$i(\mu,\gamma) = \inf_{\overline{\gamma} \in \gamma} \int_{\overline{\gamma}} d \mu $$
With this convention, we state the following corollary which will be used later.

\begin{corollary} \label{intmlmf}
Let $\gamma, \gamma' \in {\cal C}rv(S)$.
Then the following are equivalent.

\noindent
1. $\gamma \equiv_{si} \gamma'$,\\
2. $i(\lambda,\gamma) = i(\lambda,\gamma')$ for every $\lambda \in {\cal ML}(S)$, and\\
3. $i(\mu,\gamma) = i(\mu,\gamma')$ for every $\mu \in {\cal MF}(S)$.
\end{corollary}
\rightline{$\Box$}

\subsection{Representations of $\pi_{1}(S)$} \label{algprelim}

For ${\mathbb K} = {\mathbb R}$ or ${\mathbb C}$, we will consider the representation varieties
$$R_{\mathbb K}(S) = Hom(\pi_{1}(S), SL_{2}({\mathbb K})).$$
These are naturally ${\mathbb K}$-algebraic sets (see e.g. \cite{CS}).

Given an element $g \in \pi_{1}(S)$, we have the corresponding character
$$\chi_{g} : R_{\mathbb K}(S) \rightarrow {\mathbb K}$$
given by
$$\chi_{g}(\rho)=tr(\rho(g)).$$
and this is a regular function on $R_{\mathbb K}(S)$.
Since $\chi_{g}$ is invariant under conjugation (and taking inverses), it depends only on $g$ up to conjugacy (and taking inverses) in $\pi_{1}(S)$.
Since two loops represent conjugate elements in $\pi_{1}(S)$ if and only if they are {\em freely} homotopic, it follows that if $g \in \pi_{1}(S)$ is any representative of $\gamma \in {\cal C}rv(S)$ we can define $\chi_{\gamma} = \chi_{g}$.

We will also want to consider the quotient space
$$p_{\mathbb K} : R_{\mathbb K}(S) \rightarrow X_{\mathbb K}(S)$$
under the action of $GL_{2}({\mathbb K})$ by conjugation (this is almost the character variety).
The natural inclusion $SL_{2}({\mathbb R}) \subset SL_{2}({\mathbb C})$ gives us the inclusions $R_{\mathbb R}(S) \subset R_{\mathbb C}(S)$ and $X_{\mathbb R}(S) \subset X_{\mathbb C}(S)$ (this last inclusion requires a little care--this is why we quotient by the action of $GL_{2}({\mathbb K})$ instead of $SL_{2}({\mathbb K})$).

Given a point $X \in {\cal T}eich(S)$, we have the holonomy representation
$$\rho_{X}:\pi_{1}(S) \rightarrow PSL_{2}({\mathbb R})$$
which lifts to a representation (see \cite{CS})
$$\widehat{\rho}_{X}:\pi_{1}(S) \rightarrow SL_{2}({\mathbb R}).$$
The representation $\rho$ is well defined up to conjugation in $PGL_{2}({\mathbb R})$, and given $\rho$, the lifts $\widehat{\rho}$ are indexed by $H^{1}(S;{\mathbb Z}/2{\mathbb Z})$.
Choosing continuously varying lifts, we obtain a map ${\cal T}eich(S) \rightarrow X_{\mathbb R}(S)$.
Two different points $X  ,  Y \in{\cal T}eich(S)$ have non-conjugate holonomy representations so that this map is injective and we can view ${\cal T}eich(S) \subset X_{\mathbb R}(S)$.
With this identification ${\cal T}eich(S)$ sits inside $QF(S) \subset X_{\mathbb C}(S)$, the set of quasi-Fuchsian representations (discrete, faithful, geometrically finite representations without parabolics).
By Marden's Quasi-conformal Stability Theorem \cite{Mar}, ${\cal T}eich(S) \subset X_{\mathbb R}(S)$ and $QF(S) \subset X_{\mathbb C}(S)$ are open sets.

In fact, ${\cal T}eich(S)$ inherits a natural real analytic structure and $QF(S)$ inherits a natural complex analytic structure (using Fricke's coordinates, see e.g. \cite{IT} or \cite{Kquake}).
Near any $X \in {\cal T}eich(S)$ there are analytic coordinates so that ${\cal T}eich(S) \subset QF(S)$ looks like
$$U \cap {\mathbb R}^{6g-6} \subset U \cap {\mathbb C}^{6g-6}$$
where $U \subset {\mathbb C}^{6g-6}$ is an open set.

For each $\gamma \in {\cal C}rv(S)$, $\chi_{\gamma}$ is invariant under the action of $GL_{2}({\mathbb C})$, so it pushes forward to a well defined function on $X_{\mathbb C}(S)$.
The analytic structures on ${\cal T}eich(S)$ and $QF(S)$ make $\chi_{\gamma}$ into an analytic function.
The length function on ${\cal T}eich(S)$ is related to $\chi_{\gamma}$ by the formula
$$length_{X}(\gamma) = 2 cosh^{-1} \left( \frac{\chi_{\gamma}(X)}{2} \right) .$$

With this notation, we see that $\gamma \equiv_{tr} \gamma'$ means $\chi_{\gamma}^{2} = \chi_{\gamma'}^{2}$.
This {\em trace formula} immediately implies the following.

\begin{proposition} \label{tr2h}
For each $\gamma , \gamma' \in {\cal C}rv(S)$,
$$\gamma \equiv_{tr} \gamma' \Rightarrow \gamma \equiv_{h} \gamma'$$
\end{proposition}
\rightline{$\Box$}

\vspace{.3cm}

\begin{remark}
When it is convenient, and no confusion arises, we will use the same symbol to simultaneously denote an element of ${\cal C}rv(S)$, any of its representatives, the conjugacy class in $\pi_{1}(S)$ determined by it (up to taking inverses), and any representative of that conjugacy class.
\end{remark}

\section{Trace classes and hyperbolic equivalence}  \label{traces}

In \cite{H}, Horowitz studies representations of a free group $F$ into $SL_{2}(K)$, where $K$ is a commutative ring with $1$ and characteristic $0$.
In an example he constructs, for any $N > 0$, a collection of distinct elements $g_{1},...,g_{N} \in F$, such that $g_{i}$ and $g_{j}^{\pm 1}$ are non-conjugate for every $i \neq j$, and so that $\chi_{g_{i}}(\rho) = \chi_{g_{j}}(\rho)$ for each $i,j = 1,...,N$ and for all representations $\rho : F \rightarrow SL_{2}(K)$.
One of the key tools is the relation
$$tr(uv) + tr(uv^{-1}) = tr(u)tr(v)$$
for every $u,v \in SL_{2}({\mathbb K})$.

\begin{example} \cite{H}
In the rank-2 free group
$$F(2) = <a,b \, | \, - \, >$$
let $g = a^{2}b^{-1}ab$ and $h = a^{2}bab^{-1}$.
Since $g$ and $h$ are cyclically reduced, it is easy to check that $g$ is not conjugate to $h^{\pm 1}$.
Using the above trace relation one can check that for any representation
$$\rho : F(2) \rightarrow SL_{2}({\mathbb C})$$
we have
$$\chi_{g}(\rho) = \chi_{a}(\rho)[\chi_{ab}(\rho)\chi_{ab^{-1}}(\rho)-\chi_{b^{2}}(\rho)]-\chi_{a}(\rho) = \chi_{h}(\rho)$$
\end{example}

\begin{figure}[htb]
\begin{center}
\ \psfig{file=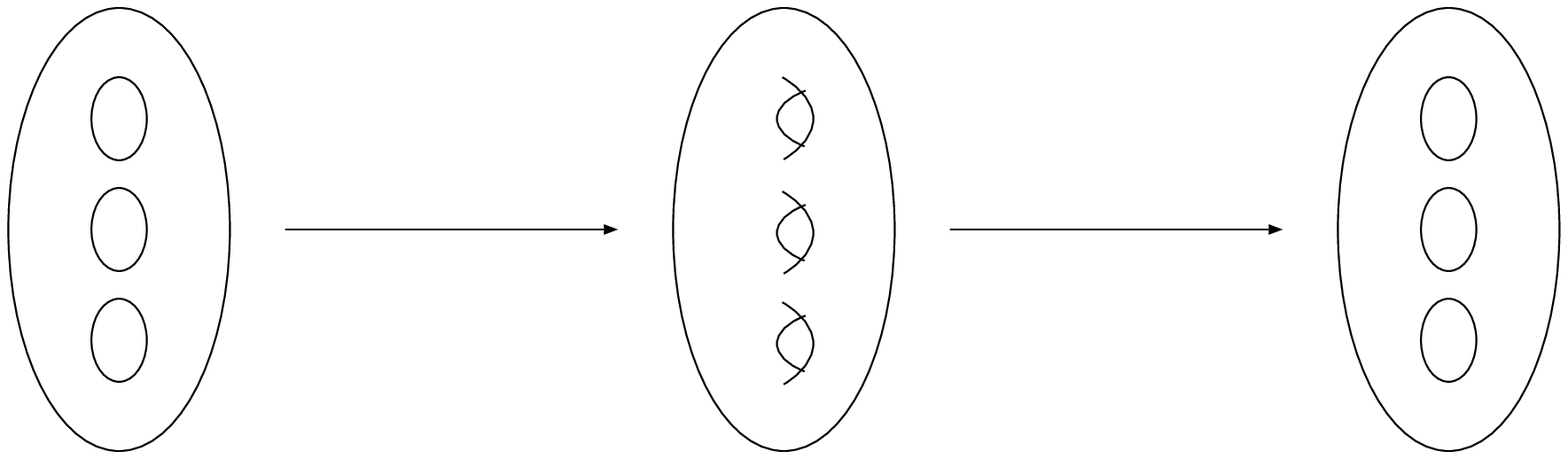,height=1truein}
\caption{$i$ and $p$.}
\label{iandp}
\end{center}
  \setlength{\unitlength}{1in}
  \begin{picture}(0,0)(0,0)
    \put(2.25,.95){$i$}
    \put(3.7,.95){$p$}
  \end{picture}
\end{figure}

Randol in \cite{Ran} uses Horowitz's examples to construct arbitrarily large $\equiv_{h}$-equivalence classes.
One way to do this is as follows.
Let $S'$ denote a sphere with $g + 1$ holes.
There are obvious maps
$$i : S' \rightarrow S \mbox{ and } p: S \rightarrow S'$$
such that $p \circ i$ is the identity (see Figure \ref{iandp}).

Let $i_{*}$ and $p_{*}$ denote the induced maps on fundamental groups.
Let $g_{1},...,g_{N} \in F(g) = \pi_{1}(S')$ be $N$ elements so that $g_{i}$ is not conjugate to $g_{j}^{\pm 1}$, for $i \neq j$, and $\chi_{g_{i}} = \chi_{g_{j}}$, for each $i,j = 1,...,N$ (as described above).
It follows that $i_{*}(g_{i})$ and $i_{*}(g_{j})^{\pm1}$ must be non-conjugate in $\pi_{1}(S)$ for each $i \neq j$, since $i_{*} \circ p_{*}$ is the identity on $F(g)$.
Therefore, $i_{*}(g_{1}),...,i_{*}(g_{N})$ represent distinct elements of ${\cal C}rv(S)$, which we refer to as $\gamma_{1},...,\gamma_{N}$ respectively, and $\gamma_{i} \equiv_{tr} \gamma_{j}$ for $i,j=1,...,N$.
By Proposition \ref{tr2h}, $\gamma_{1},...,\gamma_{N} \in {\cal C}rv(S)$ must have $\gamma_{i} \equiv_{h} \gamma_{j}$ for $i,j=1,...,N$.\\

We now begin the proof of Theorem \ref{main}.
The following proposition is the first step.
\begin{proposition} \label{tr2h2tr}
If $\gamma, \gamma' \in {\cal C}rv(S)$, then $\gamma \equiv_{tr} \gamma' \Leftrightarrow \gamma \equiv_{h} \gamma'$.
\end{proposition}
Proposition \ref{tr2h} is precisely one of the implications in the statement of this proposition.
The other implication will follow easily from the discussion in Section \ref{algprelim} and the following theorem of \cite{Rap}.

\begin{theorem} \label{rapinchuck}
(Rapinchuck et. al.) $R_{\mathbb C}(S)$ is irreducible.
\end{theorem}

\noindent
{\bf Proof of Proposition \ref{tr2h2tr}.}
It remains to show that if $\gamma \equiv_{h} \gamma'$, then $\gamma \equiv_{tr} \gamma'$.

As was noted in Section \ref{algprelim}, the functions $\chi_{\gamma}$ and $\chi_{\gamma'}$ are real and complex analytic functions on ${\cal T}eich(S)$ and $QF(S)$, respectively.
Since $length_{X}(\gamma)=length_{X}(\gamma')$ for every $X \in {\cal T}eich(S)$, the trace formula implies $\chi_{\gamma}(X) = \pm \chi_{\gamma'}(X)$ for every $X \in {\cal T}eich(S)$, and hence $\chi_{\gamma}^{2}=\chi_{\gamma'}^{2}$ on ${\cal T}eich(S)$.
Now because of the local analytic structure of the inclusion of ${\cal T}eich(S) \subset QF(S)$ described in Section \ref{algprelim}, elementary complex analysis implies $\chi_{\gamma}^{2} = \chi_{\gamma'}^{2}$ on the component of $QF(S)$ containing our embedding of ${\cal T}eich(S)$ (this then holds for any component by choosing a different embedding).

Since $QF(S) \subset X_{\mathbb C}(S)$ is open, $p_{\mathbb C}^{-1}(QF(S))$ is an open subset of $R_{\mathbb C}(S)$ (in the usual ${\mathbb C}$ topology), and thus $\chi_{\gamma}^{2} = \chi_{\gamma'}^{2}$ on this open set.
However, since $\chi_{\gamma}^{2}-\chi_{\gamma'}^{2}$ is a regular function, $\chi_{\gamma}^{2}-\chi_{\gamma'}^{2}=0$ defines an algebraic subset of $R_{\mathbb C}(S)$ containing $p_{\mathbb C}^{-1}(QF(S))$.
Some irreducible component, $Z$, of this algebraic subset contains an open subset (in the ${\mathbb C}$ topology) of $p_{\mathbb C}^{-1}(QF(S))$.
Since $Z$ and $R_{\mathbb C}(S)$ are both irreducible (the latter by Theorem \ref{rapinchuck}), they have well defined dimensions, and since $Z$ contains an open subset of $R_{\mathbb C}(S)$, these dimensions are equal.
It follows that $Z = R_{\mathbb C}(S)$ (see e.g. \cite{Shaf}), and hence $\chi_{\gamma}^{2} = \chi_{\gamma'}^{2}$ on all of $R_{\mathbb C}(S)$.\\
\rightline{$\Box$}

We state an immediate corollary of Proposition \ref{tr2h2tr}, which provides an often easily computable obstruction to $\gamma \equiv_{h} \gamma'$.

\begin{corollary} \label{h2hom}
Let $\gamma , \gamma' \in {\cal C}rv(S)$, with $\gamma \equiv_{h} \gamma'$.
Then $\gamma$ and $\gamma'$ may be oriented so that they represent the same homology class.
\end{corollary}

\noindent
{\bf Proof.}  Suppose that, regardless of the choice of orientations, $\gamma$ and $\gamma'$ represented different homology classes.
Then there exists a homomorphism $\phi : \pi_{1}(S) \rightarrow {\mathbb Z}$ such that $\phi(\gamma) \neq \pm \phi(\gamma')$.
Let $\psi : {\mathbb Z} \rightarrow SL_{2}({\mathbb C})$ be given by
$$\psi(n) = \left( \begin{array}{cc}
2^{n} & 0 \\
0 & \frac{1}{2^{n}} \\ \end{array} \right)$$
Since $tr(\psi(n)) = 2^{n} + \frac{1}{2^{n}}$, we see that $\chi^{2}_{\gamma}(\psi \circ \phi) \neq \chi_{\gamma'}^{2}(\psi \circ \phi)$, contradicting Proposition \ref{tr2h2tr}.

\rightline{$\Box$}

The following consequence of Proposition \ref{tr2h2tr} may potentially provide a topological description of hyperbolic equivalence (see Section \ref{theend}).

\begin{corollary} \label{h2maphyper}
Let $\gamma , \gamma' \in {\cal C}rv(S)$, with $\gamma \equiv_{h} \gamma'$.
Then, for every complete hyperbolic surface $M$ and (continuous) map $f : S \rightarrow M$, we have
$$length_{M}(f(\gamma)) = length_{M}(f(\gamma')).$$
\end{corollary}

\noindent
{\bf Proof.} Let $\gamma$, $\gamma'$, and $f :S \rightarrow M$ be as in the statement of the theorem.
We let $[\gamma]$ and $[\gamma']$ denote elements of $\pi_{1}(S)$ representing the conjugacy classes determined by $\gamma$ and $\gamma'$ respectively.
$f(\gamma)$ and $f(\gamma')$ determine conjugacy classes of elements of $\pi_{1}(M)$, represented by $f_{*}[\gamma]$ and $f_{*}[\gamma']$ respectively.
The hyperbolic structure determines a holonomy representation $\rho : \pi_{1}(M) \rightarrow PSL_{2}({\mathbb R})$ which lifts to a representation into $SL_{2}({\mathbb R})$ (see \cite{CS}), which we also call $\rho$.

Now, $\rho \circ f_{*}:\pi_{1}(S) \rightarrow SL_{2}({\mathbb R})$ is a homomorphism and so
$$\chi_{f_{*}[\gamma]}^{2}(\rho) = \chi_{\gamma}^{2}(\rho \circ f_{*}) = \chi_{\gamma'}^{2}(\rho \circ f_{*}) = \chi_{f_{*}[\gamma']}^{2}(\rho)$$
by Proposition \ref{tr2h2tr}.
Applying the trace formula for length we see that $f(\gamma)$ and $f(\gamma')$ have the same length in $M$.\\
\rightline{$\Box$}

\begin{remark}
Corollary \ref{h2maphyper} also holds if $M$ is a hyperbolic 3-manifold.
That is, if $\gamma \equiv_{h} \gamma'$, and $f:S \rightarrow M$ is a continuous map to a hyperbolic 3-manifold, $M$, then $f(\gamma)$ and $f(\gamma')$ have the same lengths.
\end{remark}

\section{Branched flat equivalence and simple intersection equivalence}  \label{quaddiff}

\subsection{Branched flat metrics and holomorphic quadratic differentials}

A branched flat metric on $S$ is a Euclidean cone metric such that all cone angles are of the form $k \pi$ for $k \in {\mathbb Z}$ and $k > 2$ (see e.g. \cite{CHK}), and such that the holonomy around any loop in $S$ minus the singular locus lies in $\{ 0 , \pi \}$.

A branched flat metric defines a (complete) geodesic metric on $S$, such that the metric pulled back to the universal cover has unique geodesics connecting any two points (see \cite{Ab} or \cite{BH}).
The geodesics in this metric are straight in the complement of the singular locus and make an angle no less than $\pi$ at any singular point which the geodesic hits.
Any $\gamma \in {\cal C}rv(S)$ has a length minimizing geodesic representative.

Given a complex structure $X \in {\cal T}eich(S)$, we let $Q_{X}$ denote the vector space of holomorphic quadratic differentials on $S$  (holomorphic with respect to $X$).
The union of the spaces $Q_{X}$, as $X$ ranges over ${\cal T}eich(S)$, forms a complex vector bundle over ${\cal T}eich(S)$ which we denote by $Q(S)$, and we let $Z(S)$ denote the zero section.
We claim that any $\phi \in Q(S) \setminus Z(S)$ defines a branched flat metric (see \cite{Ab}, \cite{G}, or \cite{IT}).

To see this, we suppose $\phi \in Q_{X}$ and let $z$ be a local coordinate about a point $p_{0}$ for which $\phi(p_{0}) \neq 0$.
Then $\phi = \phi(z) dz^{2}$ and in a small enough neighborhood of $p_{0}$, we can integrate to obtain a new local coordinate
$$\zeta(\omega) = \int_{z_{0}}^{\omega} \sqrt{\phi(z)} dz$$
In this new local coordinate, $\phi$ has the simple expressions $\phi = d\zeta^{2}$.
We call a local coordinate, obtained in this way, a {\em preferred coordinate}.
Any two preferred coordinates differ by rotation through $\pi$ and translation.
That is, if $\zeta_{1}$ and $\zeta_{2}$ are two preferred coordinates about a point $p$ (not a zero of $\phi$), then
$$\zeta_{1} = \pm \zeta_{2} + \omega$$
for some $\omega \in {\mathbb C}$.

It now follows that $|\phi|^{1/2}$ is the length element for a Riemannian metric at all points of $X$, except the zeros of $\phi$.
In the preferred coordinate $\zeta$, this has the simple form $|d\zeta|$, which is the standard Euclidean metric.
To understand the behavior of the metric in a neighborhood of the zero $p$ of $\phi$, we note that there exists a local coordinate $z$ about $p$ such that $\phi = z^{m}dz^{2}$.
It is now straightforward to show that the cone angles are of the required type.

Finally to see that the holonomy lies in $\{ 0, \pi \}$, we note that $\phi$ defines another type of geometric structure; a pair of transverse measured foliations.
This follows from the above description of the transition functions with respect to the preferred coordinates $\zeta = x + i y$.
The horizontal and vertical foliations with transverse measure defined by $|dy|$ and $|dx|$ respectively, are invariant by these transition functions, hence define measured foliations on $X$.
The corresponding measured foliations are called the {\em horizontal} and {\em vertical} foliations of $\phi$ and will be denoted by $({\cal F}_{h}(\phi),|Im(\sqrt{\phi})|)$ and $({\cal F}_{v}(\phi),|Re(\sqrt{\phi})|)$, respectively.
Note that the leaves of both foliations are geodesics and that they are orthogonal.
Using either one of these foliations as a ``guide'', it is clear that the holonomy lies in $\{ 0, \pi \}$.

The {\em height} of $\gamma \in {\cal C}rv(S)$ with respect to $\phi$ is defined to be
$$h_{\phi}(\gamma) = \inf_{\gamma_{0} \in \gamma} \int_{\gamma_{0}} |Im(\sqrt{\phi})|.$$
This is, of course, simply the variation of $\gamma$ with respect to the horizontal foliation.
The {\em width} is similarly defined to be the variation with respect to the vertical foliation.
It is not hard to show that if $\gamma_{0}$ is a $\phi$-geodesic representative for $\gamma$, then
$$h_{\phi}(\gamma) = \int_{\gamma_{0}} |Im(\sqrt{\phi})|.$$
Similarly for the vertical foliation.\\

Now suppose $S$ is given a branched flat structure.
This {\em defines} both a complex structure and a non-zero holomorphic quadratic differential, up to multiplication by a complex unit.
An atlas for the complex structure is described by local isometries from neighborhoods of the non-singular points to open sets in ${\mathbb E}^{2} \cong {\mathbb C}$.
In neighborhoods of the singularities, one can construct coordinate charts modeled on appropriate branched covers.

To describe the quadratic differential, note that because the holonomy lies in $\{ 0 , \pi \}$, given any non-singular point and any line in the tangent space to that point, we may parallel translate that line to every non-singular point giving a line field well defined on the complement of the singular locus.
This allows us to further refine our first type of coordinate charts above so that our line field is always sent to the line field spanned by $ \partial / \partial x$ in ${\mathbb C}$.
The transition functions for these new coordinate charts are of the form $z \mapsto \pm z + \xi$ for some $\xi \in {\mathbb C}$, and hence preserve $dz^{2}$.
Letting $\phi$ denote the pull back of $dz^{2}$ under these coordinates and extending by zero over the singular locus gives a quadratic differential, holomorphic in the complement of the singular locus.
One can check that this is in fact holomorphic everywhere.
Note that rotating the line field amounts to multiplying $\phi$ by a complex scalar with unit norm.

\subsection{Branched flat = simple intersection}

We now set out to prove the following.

\begin{proposition} \label{bf2si2bf}
If $\gamma, \gamma' \in {\cal C}rv(S)$, then $\gamma \equiv_{bf} \gamma' \Leftrightarrow \gamma \equiv_{si} \gamma'$.
\end{proposition}

We begin with the following, obvious lemma.

\begin{lemma} \label{obvious}
Let $\gamma, \gamma' \in {\cal C}rv(S)$.
If $\gamma \equiv_{si} \gamma'$ then $h_{\phi}(\gamma) = h_{\phi}(\gamma')$, for all $\phi \in Q(S) \setminus Z(S)$.
\end{lemma}

\noindent
{\bf Proof. } From the previous sub-section we know $h_{\phi}(\gamma) = i(|Re(\sqrt{\phi})|,\gamma)$.
Corollary \ref{intmlmf} implies the result.\\
\rightline{$\Box$}

\begin{remark}
It follows from the work of Hubbard and Masur in \cite{HM} that $\gamma \equiv_{si} \gamma'$ if and only if $h_{\phi}(\gamma) = h_{\phi}(\gamma')$, for all $\phi \in Q(S) \setminus Z(S)$, although we do not need this stronger version.
\end{remark}

Next, we prove the following well known lemma constituting half of Proposition \ref{bf2si2bf}.

\begin{lemma} \label{bf2si}
If $\gamma, \gamma' \in {\cal C}rv(S)$, then $\gamma \equiv_{bf} \gamma' \Rightarrow \gamma \equiv_{si} \gamma'$.
\end{lemma}

{\bf Proof.} The proof of this lemma follows from analyzing the limiting behavior of the Teichm\"uller deformation associated to an appropriate quadratic differential (see for example \cite{Kasym}).

Explicitly, we begin by fixing a curve $\alpha \in {\cal S}(S)$.
One can easily construct $\phi_{\alpha} \in Q(S) \setminus Z(S)$ such that the horizontal foliation $({\cal F}_{h}(\phi_{\alpha}),|Im(\phi_{\alpha})|)$ represents $(1,\alpha)$ under the inclusion ${\mathbb R}_{+} \times {\cal S}(S) \subset {\cal MF}(S)$.
If we remove the critical trajectories of such $\phi_{\alpha}$, the result is an open Euclidean annulus with core curve isotopic to $\alpha$, length $1$, and girth $R$, for some $R > 0$.
We deform such a $\phi_{\alpha}$ to obtain a family $\{ \phi_{\alpha,K} \}_{K \in (0,1]} \subset Q(S) \setminus Z(S)$, by compressing the horizontal direction by a factor of $K$.
More precisely, let $e_{1},e_{2}$ be an orthonormal basis for the tangent space at a non-singular point $x \in S$, with $e_{1}$ tangent to the horizontal foliation (and hence $e_{2}$ tangent to the vertical). 
The original metric at $x$ is thus given by $\delta_{ij}$, and the deformed metric is given by $g_{12} = g_{21} = 0$, $g_{11} = K^{2}$, and $g_{22} = 1$.
For every $K \in (0,1]$, removing the critical trajectories of $\phi_{\alpha,K}$ we obtain an open Euclidean annulus with length $1$ and girth $KR$, so that $({\cal F}_{h}(\phi_{\alpha,K}),|Re(\sqrt{\phi_{\alpha,K}})|)$ also represents $(1,\alpha)$.

We now let $\gamma \in {\cal C}rv(S)$ denote any element, and claim that
$$\lim_{K \rightarrow 0} length_{\phi_{\alpha,K}}(\gamma) = h_{\phi_{\alpha,K}}(\gamma) = i(\alpha,\gamma).$$
The second equality holds by construction.
To see the first, let $\gamma_{0}$ denote a $\phi_{\alpha}$-geodesic representative for $\gamma$.
Then $\gamma_{0}$ consists of $i(\alpha,\gamma)$ geodesic arcs traversing the length of the annulus and some number of geodesic arcs, each of which is contained in the set of critical trajectories.
It is not hard to see that $\gamma_{0}$ remains geodesic after the deformation, that is, it is also a $\phi_{\alpha,K}$-geodesic.
As $K \rightarrow 0$, the length of the arcs contained in the critical trajectories converges to $0$, and the length of the others converges to $1$.

To complete the proof, let $\gamma, \gamma' \in {\cal C}rv(S)$ be such that $\gamma \equiv_{bf} \gamma'$, and let $\phi_{\alpha,K}$ be as above.
Then
$$i(\alpha,\gamma) = \lim_{K \rightarrow 0} length_{\phi_{\alpha,K}}(\gamma)
=\lim_{K \rightarrow 0} length_{\phi_{\alpha,K}}(\gamma') = i(\alpha,\gamma').$$
Since $\alpha$ was arbitrary, $\gamma \equiv_{si} \gamma'$.\\
\rightline{$\Box$}

To prove the other implication of Proposition \ref{bf2si2bf} we will need the following

\begin{definition}
Given a finite list of points $Z = (z_{1},...,z_{k} ) \in {\mathbb C}^{k}$, define the {\em height} of $Z$ to be
$$h(Z) = \sum_{j=1}^{k} |Im (z_{j})|.$$
\end{definition}

\begin{remark}
We are considering lists of numbers instead of sets to allow the possibility of a number showing up with multiplicity greater than one.
\end{remark}

\begin{lemma} \label{pointheight}
Let $Z = (z_{1},...,z_{k}) \in {\mathbb C}^{k}$ be a finite list of points.
Then $H_{Z}(\theta) = h(e^{i\theta} Z)$ is continuous and piecewise smooth.
If $\theta_{1},...,\theta_{n}$ denotes the points in $(0,\pi]$ where $H_{Z}'(\theta)$ does not exists, then
$$ \sum_{j=1}^{k} |z_{j}| = \frac{1}{2} \sum_{l = 1}^{n} \left( \lim_{\theta \rightarrow \theta^{+}_{l}} H_{Z}'(\theta) - \lim_{\theta \rightarrow \theta^{-}_{l}} H_{Z}'(\theta) \right) .$$
\end{lemma}

\noindent
{\bf Proof.} $H_{Z}$ is clearly continuous.

To prove the remainder of this lemma, it will be convenient to define the function $f : {\mathbb C} \rightarrow {\mathbb C}$ by
$$ f(z) = \left\{ \begin{array}{cc}
z & \mbox{ if } Im(z) > 0 \mbox{ or } Im(z) = 0 \mbox{ and } Re(z) \geq 0 \\
-z & \mbox { otherwise }\\ \end{array} \right\}.$$
Then
$$H_{Z}(\theta) = \sum_{j=1}^{k} Im ( f(e^{i \theta}z_{j}) ).$$

Note that for any $z \in {\mathbb C}^{*}$, the function $h(e^{i \theta}z) = Im( f (e^{i \theta}z) )$ is smooth (as a function of $\theta$) exactly when $ e^{i \theta}z \not\in {\mathbb R}$, i.e. when $arg(z) + \theta \not\in {\mathbb Z}\pi$.
Moreover, away from such values of $\theta$, the derivative is $\frac{d}{d \theta} h(e^{i \theta}z) = \frac{d}{d \theta} Im( f(e^{i \theta} z) ) =  Re( f(e^{i \theta} z) )$.
If we let $\theta_{z}$ be any value of ${\mathbb Z} \pi - arg(z)$, then
$$\lim_{\theta \rightarrow \theta^{+}_{z}} f(e^{i \theta}z) = |z|
\mbox{ and }
\lim_{\theta \rightarrow \theta^{-}_{z}} f(e^{i \theta}z) = - |z|$$
As $\theta \rightarrow \theta_{z}$ from the right we have $Im( f( e^{i \theta} z) ) \rightarrow 0$ so that
$$ \lim_{\theta \rightarrow \theta^{+}_{z}} \frac{d}{d \theta} ( Im( f( e^{i \theta}z) ) )
= \lim_{\theta \rightarrow \theta^{+}_{z}} Re( f( e^{i \theta}z) )
= \lim_{\theta \rightarrow \theta^{+}_{z}} f(e^{i \theta}z)
= |z|$$
Putting this together with the analogous result when $\theta \rightarrow \theta_{z}$ from the left, we obtain
\begin{equation} \label{onlyequation}
\lim_{\theta \rightarrow \theta^{+}_{z}} \frac{d}{d \theta} ( Im( f( e^{i \theta}z) ) ) -  \lim_{\theta \rightarrow \theta^{-}_{z}} \frac{d}{d \theta} ( Im( f( e^{i \theta}z) ) ) = |z| - (- |z|) = 2 |z|
\end{equation}

We now return to the proof.
We may clearly assume that no coordinate of $Z$ is $0$.
Also, note that we may replace $Z$ by $(f(z_{1}),...,f(z_{k}))$ without loss of generality, so that all coordinates of $Z$ lie in the open upper half plane union the positive real axis.
Further, if $z_{i}$ is a (positive) real multiple of $z_{j}$ ($j \neq i$), then we may replace the two coordinates, $z_{i}$ and $z_{j}$, by the single coordinate $z_{i}+z_{j}$.
We therefore assume that no two coordinates lie on the same line through $0$ in ${\mathbb C}$.

It follows that the arguments $Arg(z_{1}),...,Arg(z_{k})$, taken in $[0,\pi)$ are distinct.
$H_{Z}$ is now seen to be smooth except when $\theta = m \pi - Arg(z_{j})$ and $m \in {\mathbb Z}$.
Letting $\theta_{j} = \pi - Arg(z_{j})$ for each $j = 1,...,k$, $\{ \theta_{1},...,\theta_{k} \}$ is precisely the set of points in $(0,\pi]$ where $H_{Z}$ is not smooth.

We have
$$\sum_{j=1}^{k} \left( \lim_{\theta \rightarrow \theta^{+}_{j}} H_{Z}'(\theta) - \lim_{\theta \rightarrow \theta^{-}_{j}} H_{Z}'(\theta) \right)
= \sum_{j=1}^{k} \sum_{l=1}^{k} \left( \lim_{\theta \rightarrow \theta^{+}_{j}}(\frac{d}{d \theta} Im( f( e^{i \theta} z_{l} ) ) ) - \lim_{\theta \rightarrow \theta^{-}_{j}}(\frac{d}{d \theta} Im( f( e^{i \theta} z_{l} ) )) \right)$$
$$= \sum_{j=1}^{k} 2 | z_{j} |$$
where the last equality follows from equation (\ref{onlyequation}) when $j = l$ and from the fact that $Im( f(e^{i \theta} z_{l} ) )$ is smooth at $\theta_{j}$ whenever $j \neq l$.
This completes the proof of the lemma.\\
\rightline{$\Box$}

The following lemma, which seems interesting in its own right, will easily complete the proof of Proposition \ref{bf2si2bf}.

\begin{lemma} \label{lengthfromheight}
Let $\gamma \in {\cal C}rv(S)$ and $\phi \in Q(S) \setminus Z(S)$.
Then $H_{\gamma}(\theta) = h_{e^{2 i \theta} \phi}(\gamma)$ is continuous and piecewise smooth.
If $\theta_{1},...,\theta_{n}$ denotes the set of points in $(0,\pi]$ for which $H_{\gamma}'$ does not exist, then
$$length_{\phi}(\gamma) = \frac{1}{2} \sum_{l=1}^{n} \left( \lim_{\theta \rightarrow \theta^{+}_{l}} H_{\gamma}'(\theta) - \lim_{\theta \rightarrow \theta^{-}_{l}} H_{\gamma}'(\theta) \right) .$$
\end{lemma}
 
\noindent
{\bf Proof.}
Let $\gamma_{0}$ denote a $\phi$-geodesic representative for $\gamma$.
This is also an $e^{2i \theta} \phi$-geodesic, since the metric defined by $e^{2i \theta} \phi$ is independent of $\theta$.
Let $\gamma_{0,1}, ..., \gamma_{0,k}$ denote the straight arcs (with endpoints on the zeros of $\phi$) which make up $\gamma_{0}$.
For every $\theta \in {\mathbb R}$ and $j \in \{1,...,k \}$, we can develop $\gamma_{0,j}$ to a straight arc $\gamma_{0,j,\theta} \subset {\mathbb C}$ using preferred coordinates defined by $e^{2i \theta} \phi$.

If $\gamma_{0,j,\theta}$ and $\gamma_{0,j,\theta}'$ are any two developing images of $\gamma_{0,j}$ by $e^{2i \theta} \phi$, then there is $\omega \in {\mathbb C}$, such that $\gamma_{0,j,\theta}' = \psi(\gamma_{0,j,\theta})$, where $\psi(\zeta) = \pm \zeta + \omega$.
We may therefore assume that each $\gamma_{0,j,\theta}$ has one endpoint at $0$.
Denote the other endpoint of $\gamma_{0,j,\theta} = z_{j,\theta}$.

Set $Z_{\theta} = (z_{1,\theta}, ... , z_{k,\theta})$ and note that $h(Z_{\theta}) = h_{e^{2i \theta}\phi}(\gamma)$ since the horizontal foliation of $e^{2i \theta}\phi$ is obtained by pulling back the horizontal foliation of ${\mathbb C}$ by a preferred coordinate, and since geodesics realize the height.

If we let $\zeta$ denote a preferred coordinate for $\phi$, then $e^{i \theta} \zeta$ is a preferred coordinate for $e^{2i \theta}\phi$.
Therefore, we can take $Z_{\theta} = (e^{i \theta}z_{1,0},...,e^{i \theta}z_{k,0}) = e^{i \theta} Z_{0}$.

So, $H_{\gamma}(\theta) = h(Z_{\theta}) = h(e^{i \theta} Z_{0})=H_{Z_{0}}(\theta)$.
Since
$$length_{\phi}(\gamma) = length_{\phi}(\gamma_{0}) = \sum_{j=1}^{k} length_{\phi} (\gamma_{0,j}) = \sum_{j=1}^{k} length_{\mathbb C} (\gamma_{0,j,0}) = \sum_{j=1}^{k} |z_{j,0}|$$
Lemma \ref{pointheight} applied to $Z_{0}$ completes the proof.

\rightline{$\Box$}

\noindent
{\bf Proof of Proposition \ref{bf2si2bf}. }
By Lemma \ref{bf2si}, we need only prove that if $\gamma, \gamma' \in {\cal C}rv(S)$ and $\gamma \equiv_{si} \gamma'$, then $\gamma \equiv_{bf} \gamma'$.

Fix $\phi \in Q(S) \setminus Z(S)$ arbitrarily.
Let $H_{\gamma}(\theta)$ and $H_{\gamma'}(\theta)$ be as in the statement of Lemma \ref{lengthfromheight}.
By Lemma \ref{obvious} , $H_{\gamma}(\theta) = H_{\gamma'}(\theta)$ for every $\theta \in {\mathbb R}$, and so applying Lemma \ref{lengthfromheight}, we obtain $length_{\phi}(\gamma) = length_{\phi}(\gamma')$.
Since $\phi$ was arbitrary, we are done.

\rightline{$\Box$}

\begin{remark}
As a consequence of this, Corollary \ref{h2si} below, and the existence of arbitrarily large hyperbolic classes in ${\cal C}rv(S)$, we see that there are arbitrarily large branched flat classes in ${\cal C}rv(S)$.
\end{remark}

\section{Geodesic currents} \label{geodesiccurrents}

In this section, we briefly discuss Thurston's compactification of Teichm\"uller space in the context of Bonahon's work on geodesic currents \cite{BTeich} as it applies to our situation.

Following \cite{BTeich}, we define a {\em geodesic current} as follows.
We have on $S$ a fixed hyperbolic metric (cf. Section \ref{teichprelim}), which provides us a hyperbolic metric on the universal cover $p: \widetilde{S} \rightarrow S$ (thus making it isometric to ${\mathbb H}^{2}$).
We denote the circle at infinity by $S_{\infty}^{1}$ and the space of geodesics on $\widetilde{S}$ by $G(\widetilde{S})$.
The endpoints of a geodesic naturally provide an identification of $G(\widetilde{S})$ with $(S_{\infty}^{1} \times S_{\infty}^{1} \setminus \Delta) / \sim$ where $\Delta \subset S_{\infty}^{1} \times S_{\infty}^{1}$ is the diagonal and $(x,y) \sim (y,x)$.
The fundamental group $\pi_{1}(S)$ acts isometrically on $\widetilde{S}$, and hence it also acts on $G(\widetilde{S})$.
We define a {\em geodesic current} on $S$ to be a $\pi_{1}(S)$-invariant positive Borel measure on $G(\widetilde{S})$ (measuring compact sets finitely).
We denote the space of geodesic currents on $S$ by ${\cal C}urr(S)$, topologized with the weak* topology.

We can naturally embed ${\mathbb R}_{+} \times {\cal C}rv(S)$ into ${\cal C}urr(S)$ as follows.
The metric on $S$ allows us to view $\gamma \in {\cal C}rv(S)$ as a geodesic so that $p^{-1}(\gamma)$ is a $\pi_{1}(S)$-invariant discrete subset $X_{\gamma} \subset G(\widetilde{S})$.
We define the corresponding geodesic current (denoted simply by $t \cdot \gamma$) to be $t$ times the atomic measure corresponding to $X_{\gamma}$ (so, $t \cdot \gamma(E) = t \cdot | E \cap X_{\gamma} |$, for any Borel set $E$).

We may also embed ${\cal ML}(S)$ into ${\cal C}urr(S)$.
The measured lamination $(\Lambda,\lambda)$ pulls back to a measured lamination $(\widetilde{\Lambda},\widetilde{\lambda})$ on $\widetilde{S}$.
We view $\widetilde{\Lambda}$ as a (closed) $\pi_{1}(S)$-invariant subset of $G(\widetilde{S})$, and define a measure on $\widetilde{\Lambda}$ as follows.
A geodesic arc $\alpha$ in $\widetilde{S}$ defines a subset $E_{\alpha} \subset \widetilde{\Lambda}$ consisting of those geodesics transversely intersecting $\alpha$.
We define the measure of $E_{\alpha}$ to be $\widetilde{\lambda}(E_{\alpha})$.
The sets of the form $E_{\alpha}$ are enough to well define a measure on $\widetilde{\Lambda}$ which we extend to all of $G(\widetilde{S})$ and denote simply by $\lambda \in {\cal C}urr(S)$.

Next, we briefly describe the (proper) embedding of ${\cal T}eich(S)$ into ${\cal C}urr(S)$ as defined in \cite{BTeich}.
The {\em Liouville measure}, $L$, is the $Isom({\mathbb H}^{2})$-invariant positive Borel measure on $G({\mathbb H}^{2})$ (the space of geodesics on ${\mathbb H}^{2}$) with the following defining property.
If $\alpha$ and $\beta$ are two disjoint arcs on the circle at infinity with endpoints $\alpha_{0}$, $\alpha_{1}$ and $\beta_{0}$, $\beta_{1}$ respectively (in either the disk or upper half plane model) and $E_{\alpha,\beta}$ is the set of geodesics with one endpoint in $\alpha$ and the other in $\beta$, then
$$L(E_{\alpha,\beta}) = \left| log \left( \left| \frac{(\alpha_{0}-\beta_{0})(\alpha_{1}-\beta_{1})}{(\alpha_{0}-\beta_{1})(\alpha_{1}-\beta_{0})} \right| \right) \right|$$
Now, given a hyperbolic structure $X$ on $S$, we view $X$ as a hyperbolic surface along with a diffeomorphism $f:S \rightarrow X$.
The map $f$ lifts $\pi_{1}(S)$-equivariantly to a diffeomorphism of the universal covers
$$\widetilde{f}: \widetilde{S} \rightarrow {\mathbb H}^{2}.$$
This map extends (by a $\pi_{1}(S)$-equivariant homeomorphism) to a map on circles at infinity, and hence induces a homeomorphism
$$\overline{f}_{X} : G(\widetilde{S}) \rightarrow G({\mathbb H}^{2})$$
Pulling back the Liouville measure via $\overline{f}_{X}$ describes a Borel measure on $G(\widetilde{S})$ which is invariant under $\pi_{1}(S)$.
That is, $L_{X} = \overline{f}_{X}^{*}(L)$ is a geodesic current.
The map
$${\cal T}eich(S) \rightarrow {\cal C}urr(S)$$
defined by $X \mapsto L_{X}$ is a proper embedding \cite{BTeich}.

For the remainder of this section, we will identify ${\mathbb R}_{+} \times {\cal C}rv(S)$, ${\cal ML}(S)$, and ${\cal T}eich(S)$ with their images in ${\cal C}urr(S)$.\\

In \cite{Bbouts}, a bilinear function is defined
$$I: {\cal C}urr(S) \times {\cal C}urr(S) \rightarrow {\mathbb R}$$
This function enjoys several nice properties (for the definition of $I$ and a proof of the following see \cite{Bbouts} and \cite{BTeich} (there, $I$ is called $i$)).

\begin{theorem} (Bonahon) \label{bonahonI} \

\noindent
1. The function $I$ is continuous.

\noindent
2. Given $t \cdot \gamma \, , \, t' \cdot \gamma' \in {\mathbb R}_{+} \times {\cal C}rv(S)$, $\lambda \in {\cal ML}(S)$, and $X \in {\cal T}eich(S)$, we have
\begin{itemize}
\item $I(t \cdot \gamma, t' \cdot \gamma') = tt' \cdot i(\gamma,\gamma')$
\item $I(\lambda,t \cdot \gamma) = t \cdot i(\lambda,\gamma)$
\item $I(L_{X},t \cdot \gamma) = t \cdot length_{X}(\gamma)$
\item $I(L_{X},L_{X}) = \pi^{2}|\chi(S)|$
\end{itemize}

\noindent
3. If $\mu \in {\cal C}urr(S) \setminus \{ 0 \}$ satisfies $I(\mu,\mu)=0$, then $\mu \in {\cal ML}(S)$.
\end{theorem}

Using this function and this embedding of ${\cal T}eich(S)$, one can recover Thurston's compactification (see \cite{Tdiff} and \cite{FLP} for the original definition and proof).
Briefly, one first notes that by the final item in Theorem \ref{bonahonI}(2), the inclusion of ${\cal T}eich(S)$ into ${\cal C}urr(S)$ remains embedded after passing to the projectivization of ${\cal C}urr(S)$, denoted ${\cal PC}urr(S) = ({\cal C}urr(S) \setminus \{ 0 \}) / {\mathbb R}_{+}$.
The space ${\cal PC}urr(S)$ is compact, and therefore the closure of the image of ${\cal T}eich(S)$ in ${\cal PC}urr(S)$ is compact (see \cite{BTeich}).

Now suppose $\{ X_{i} \}$ is a sequence in ${\cal T}eich(S)$ which diverges.
We wish to identify the limit point in ${\cal PC}urr(S)$ as an element of ${\cal PML}(S)$, the projectivization of ${\cal ML}(S)$.
By passing to a subsequence (using compactness), we may assume that the image of the $\{ X_{i} \}$ in ${\cal PC}urr(S)$ converges.
Therefore, there exists $t_{i} \in {\mathbb R}^{+}$ so that $t_{i} \cdot X_{i} \rightarrow \lambda \in {\cal C}urr(S) \setminus \{ 0 \}$.
Since the embedding of ${\cal T}eich(S)$ into ${\cal C}urr(S)$ is proper (or by an easy geometric argument), it must be that $t_{i} \rightarrow 0$.
So by Theorem \ref{bonahonI}
$$t_{i}^{2} \pi^{2} |\chi(S)| = t_{i}^{2} I(X_{i},X_{i}) =  I(t_{i}X_{i},t_{i}X_{i}) \rightarrow I(\lambda,\lambda)$$
as $i \rightarrow \infty$.
Since the first term is obviously approaching $0$, we see that $I(\lambda,\lambda) = 0$ and therefore $\lambda \in {\cal ML}(S)$.

As in \cite{BTeich}, one can check that this is exactly Thurston's compactification.

\begin{theorem} (Thurston) \label{thurstoncompact}
${\cal T}eich(S)$ is compactified by ${\cal PML}(S)$.
In this compactification, which we denote by $\overline{{\cal T}eich(S)}$, a sequence $\{ X_{i} \} \subset {\cal T}eich(S)$ converges to $[\lambda] \in {\cal PML}(S)$ if and only if for every pair $\gamma, \gamma' \in {\cal C}rv(S)$ we have
$$\lim_{i \rightarrow \infty} \frac{length_{X_{i}}(\gamma)}{length_{X_{i}}(\gamma')} = \frac{i(\lambda,\gamma)}{i(\lambda,\gamma')}$$
provided $i(\lambda,\gamma') \neq 0$.
Moreover, ${\cal T}eich(S)$ is a dense open set in $\overline{{\cal T}eich(S)}$.
\end{theorem}
\rightline{$\Box$}

\begin{remark}
In the usual statement of this theorem $\gamma$ and $\gamma'$ are taken to lie in ${\cal S}(S)$, but the work in \cite{BTeich} described above easily implies this version.
The usual statement also includes the additional information that $\overline{{\cal T}eich(S)} \cong \overline{B}^{6g-6}$, the closed ball in ${\mathbb R}^{6g-6}$, and that ${\cal T}eich(S)$ is the interior (see \cite{Tdiff} or \cite{FLP}).
\end{remark}

The following is now immediate from this theorem.

\begin{corollary} \label{h2si}
Let $\gamma,\gamma' \in {\cal C}rv(S)$.
Then $\gamma \equiv_{h} \gamma' \Rightarrow \gamma \equiv_{si} \gamma'$.
\end{corollary}

\noindent
{\bf Proof. } Let $\alpha \in {\cal S}(S)$ be an arbitrary simple closed curve, and let us denote the image of $(1,\alpha)$ in ${\cal ML}(S)$ by $\alpha$.
By Theorem \ref{thurstoncompact}, there is a sequence $X_{i} \in {\cal T}eich(S)$ such that $X_{i} \rightarrow [\alpha]$ (it is not difficult to explicitly construct such a sequence).
Since $\gamma \equiv_{h} \gamma'$, we have
$$1 = \lim_{i \rightarrow \infty} \frac{length_{X_{i}}(\gamma)}{length_{X_{i}}(\gamma')} = \frac{i(\alpha,\gamma)}{i(\alpha,\gamma')}$$
which completes the proof.\\
\rightline{$\Box$}\\

Now, let $\gamma,\gamma' \in {\cal C}rv(S)$, and consider the following function
$$\Phi_{\gamma,\gamma'} : {\cal T}eich(S) \rightarrow {\mathbb R}$$
given by
$$\Phi_{\gamma,\gamma'}(X) = \frac{length_{X}(\gamma)}{length_{X}(\gamma')}$$

From Section \ref{algprelim}, we see that $\Phi_{\gamma,\gamma'}$ is not only continuous on ${\cal T}eich(S)$, but it is real analytic.
This function does not necessarily extend continuously over $\overline{{\cal T}eich(S)}$ in general (even when $\gamma \equiv_{si} \gamma'$).

However, if we assume that both $\gamma$ and $\gamma'$ (individually) fill $S$, so that for every $\lambda \in {\cal ML}(S)$, $i(\lambda,\gamma) \neq 0 \neq i(\lambda,\gamma')$, then there is an extension (denoted $\overline{\Phi}_{\gamma,\gamma'}$) which is continuous by Theorem \ref{mlextend} and Theorem \ref{thurstoncompact}

Moreover, if in addition we assume that $\gamma \equiv_{si} \gamma'$, then $\overline{\Phi}_{\gamma,\gamma'} = 1$ on ${\cal PML}(S)$.
In particular, $\Phi_{\gamma,\gamma'}$ is asymptotic to $1$ in every direction, and hence $\Phi_{\gamma,\gamma'}$ extends by $1$ over {\em any} compactification of ${\cal T}eich(S)$.

Despite this behavior, in Section \ref{simpleintersection} we construct pairs of curves $\gamma$ and $\gamma'$ so that $\gamma$ and $\gamma'$ both fill the surface and $\gamma \equiv_{si} \gamma'$, yet $\gamma \not \equiv_{h} \gamma'$.

\section{Simple intersection equivalence} \label{simpleintersection}

In this section, we construct a pair of closed curves which both (individually) fill $S$ and are simple intersection equivalent, yet are not hyperbolically equivalent.
The first difficulty that one encounters in constructing such a pair of curves is verifying that they {\em are} simple intersection equivalent.
We discuss this difficulty at some length as it strongly contrasts the behavior of simple curves.\\

\subsection{Distinguishing curves}

The work in this sub-section is not necessary for the proof of the main theorem, but is included because of the obvious relevance to constructing the required counterexamples.

It is well known (see e.g. \cite{PH}) that for any surface $S$, there exists a finite set $\{ \alpha_{1},...,\alpha_{n} \} \subset {\cal S}(S)$ such that if $\gamma, \gamma' \in {\cal S}(S)$, and
$$i(\alpha_{j},\gamma) = i(\alpha_{j},\gamma')$$
for each $j = 1,..., n$, then $\gamma = \gamma'$.

Of course, when $\gamma, \gamma' \in {\cal C}rv(S)$ (and $S$ has genus at least $2$) this fails by Corollary \ref{h2si} and the existence of $\gamma \equiv_{h} \gamma'$ with $\gamma \neq \gamma'$ (Section \ref{traces}).
In fact the situation is worse, as the following indicates.

\begin{theorem} \label{thatswhy}
Given any proper compact subset $K \subset {\cal PML}(S)$ there exists two curves $\gamma,\gamma' \in {\cal C}rv(S)$ and $\alpha \in {\cal S}(S)$ such that
$$i(\lambda,\gamma) = i(\lambda,\gamma')$$
for every $\lambda$ with $[\lambda] \in K$, yet $i(\alpha,\gamma) \neq i(\alpha,\gamma')$
\end{theorem}

In particular, this theorem says that there is {\em no} finite list of simple closed curves such that for any $\gamma, \gamma' \in {\cal C}rv(S)$, the intersection number with each curve in this list can be used to decide whether or not $\gamma$ and $\gamma'$ are simple intersection equivalent.\\

\noindent
{\bf Proof.}
Since $K$ is a proper compact subset of ${\cal PML}(S)$ there exists a lamination $\lambda \in {\cal ML}(S)$ such that $[\lambda] \not \in K$ and so that some complementary region of $\lambda$ is an ideal triangle.
It follows from \cite{Tnotes} (in particular \S 9.7) that $\lambda$ may be chosen so that there are train tracks $\tau$ and $\tau_{K}$ with $\tau$ carrying $\lambda$, $\tau_{K}$ carrying every element of $K$, and so that $\tau$ and $\tau_{K}$ meet efficiently, in the sense that they meet transversely and so that any trainpaths of $\tau$ and $\tau_{K}$ lift to paths in the universal cover which meet transversely in at most one point.
This efficiency guarantees that if $\gamma$ is carried by $\tau$ and $\gamma_{K}$ is carried by $\tau_{K}$, then $i(\gamma,\gamma_{K})$ is determined by the respective weights they define on the branches of $\tau$ and $\tau_{K}$.
In fact, it is not hard to see that this holds even if we allow $\gamma$ to be a non-simple closed curve carried by $\tau$ (the definition of a train track carrying a closed curve is essentially the same as that given for a lamination in \cite{Tnotes}).

\begin{figure}[htb]
\begin{center}
\ \psfig{file=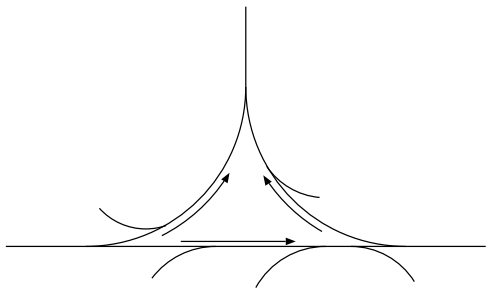,height=1.5truein}
\caption{Incompatible.}
\label{incompatible}
\end{center}
\end{figure}

Now let $\gamma$ be an oriented simple closed curve fully carried by $\tau$ (i.e. $\gamma$ defines positive weights on all branches of $\tau$).
Note that $\tau$ must have a complementary region which is a triangle.
It follows that there must be some branch which contains strands of $\gamma$ with disagreeing orientations.
To see this, note that if this weren't the case, then the boundary of the triangular regions could be oriented to agree with the orientation of $\gamma$.
Thus two sides of the triangle would have to meet at a vertex with one orientation heading into the vertex and the other heading out (see Figure \ref{incompatible}).
A branch on the opposite side of this switch would have to contain strands with both orientations (because the switches are trivalent).

\begin{figure}[htb]
\begin{center}
\ \psfig{file=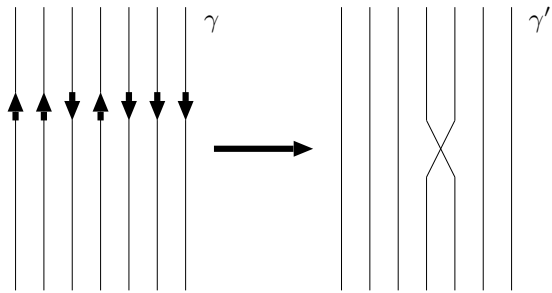,height=1.5truein}
\caption{Exchange.}
\label{exchange}
\end{center}
  \setlength{\unitlength}{1in}
  \begin{picture}(0,0)(0,0)
    \put(2.73,2){$\gamma$}
    \put(4.43,2){$\gamma'$}
  \end{picture}
\end{figure}

Next choose two strands in such a branch with opposite orientation which are adjacent (here we are thinking of $\gamma$ as embedded in a small ``$I$-bundle'' neighborhood of $\tau$).
Cut these strands and re-glue them with a cross (see Figure \ref{exchange}).
We claim that the resulting immersed $1$-manifold is actually a curve, $\gamma' \in {\cal C}rv(S)$ and it is not in ${\cal S}(S)$ (i.e. it is not homotopic to a {\em simple} curve).

Assuming this claim for the moment, we see that $\gamma'$ is also carried by $\tau$, and that $\gamma$ and $\gamma'$ define the same weights on the branches of $\tau$.
It follows that for any measured lamination $\mu$ carried by $\tau_{K}$ (in particular, any $\mu \in K$), that $i(\mu,\gamma) = i(\mu,\gamma')$.
Because $\gamma \in {\cal S}(S)$ and $\gamma'$ is not, one can easily find $\alpha \in {\cal S}(S)$ such that $i(\alpha,\gamma) \neq i(\alpha,\gamma')$.

To complete the proof of the proposition, we must prove:\\\\
i. $\gamma'$ has one component (as an immersed $1$-manifold), and \\\\
ii. the single self-intersection point cannot be removed.\\

\noindent
{\bf proof of i. } Let $\gamma_{0} = \gamma_{1} \amalg \gamma_{2}$ denote the $2$ {\em oriented} arcs which remain after cutting $\gamma$.
Re-gluing will result in a disconnected $1$-manifold if the two ends of $\gamma_{1}$ are glued together (and hence also the two ends of $\gamma_{2}$).
Since we have chosen a pair of strands which are oppositely oriented, the negative boundary components of $\gamma_{0}$ are glued together and hence $\gamma_{1}$ is glued to $\gamma_{2}$.
Therefore, $\gamma'$ is connected and (i) holds.\\

\noindent
{\bf proof of ii. } Again, we think of $\gamma$ as embedded in an $I$-bundle neighborhood of $\tau$.
By construction, $\gamma'$ has exactly one (transverse) self intersection point.
An exercise in the topology of surfaces shows that a curve with one transverse self intersection point is homotopic to a simple closed curve if and only if one of the complementary regions of the curve is a disk (the complementary regions being the components of the path metric completion of the surface minus the curve).
The components of the preimage of $\gamma'$ in $\widetilde{S}$ (i.e. the complete lifts of $\gamma'$) are easily seen to be embedded (a component of the preimage is a trainpath and so the intersection with the $I$-bundle neighborhood of a branch of $p^{-1}(\tau)$ is a single arc).
A complementary region which is a disk would lift to the universal cover and no component of the preimage could be embedded.
Therefore, no complementary region is a disk and hence $\gamma'$ is not homotopic to a simple closed curve.\\
\rightline{$\Box$}

\subsection{Pairs of pants}

There is a situation where we can easily decide whether or not a pair of curves are simple intersection equivalent without knowing that they are hyperbolically equivalent.
We will use this in our construction in the next subsection.

\begin{figure}[htb]
\begin{center}
\ \psfig{file=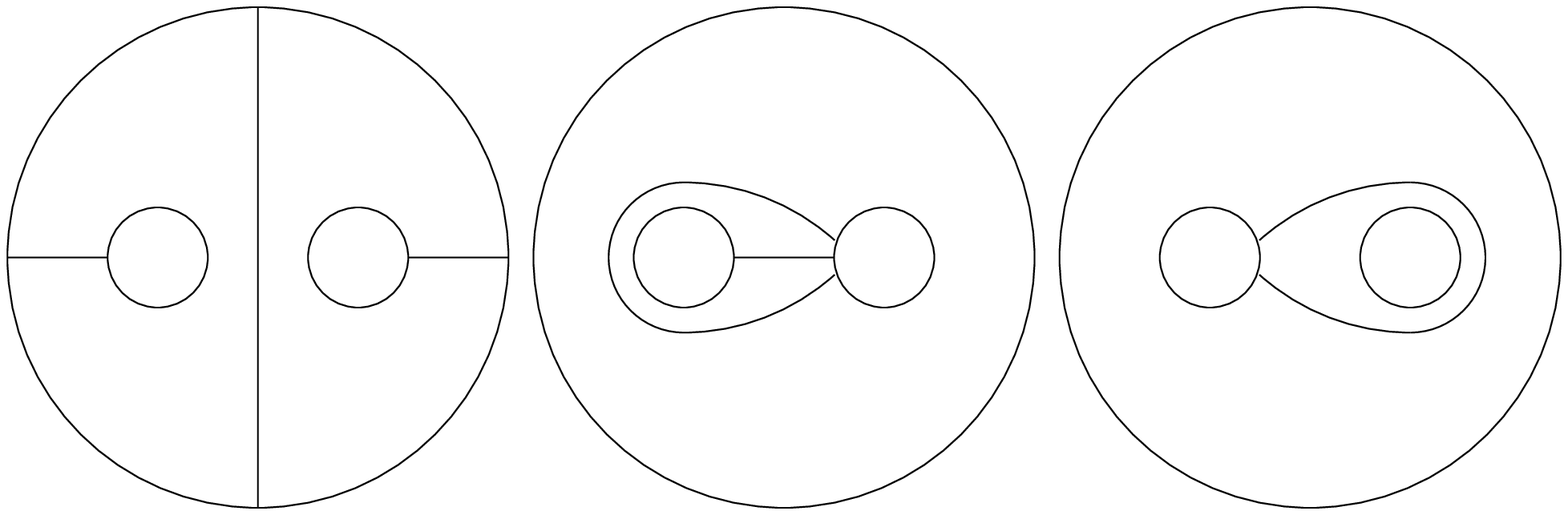,height=1.3truein}
\caption{Arcs in Pants.}
\label{arcsinpants}
\end{center}
  \setlength{\unitlength}{1in}
  \begin{picture}(0,0)(0,0)
    \put(1.15,1.35){$l_{1}$}
    \put(1.7,1.7){$w_{1}$}
    \put(2.1,1.35){$l_{2}$}
    \put(2.9,1.32){$l_{3}$}
    \put(2.45,1.1){$w_{2}$}
    \put(4.8,1.4){$w_{3}$}
  \end{picture}
\end{figure}

Let $P$ denote a pair of pants, i.e. a sphere with $3$ holes.
There are exactly $6$ isotopy classes of essential properly embedded arcs which we call $l_{1},l_{2},l_{3}$ and $w_{1},w_{2}, w_{3}$ (see \cite{PH}).
Representatives of each isotopy class are shown in Figure \ref{arcsinpants}.
Note that for each $i=2,3$ there is a homeomorphism of $P$ taking $l_{1}$ to $l_{i}$.
Similarly, for each $i=2,3$ there is a homeomorphism of $P$ taking $w_{1}$ to $w_{i}$.

\begin{figure}[htb]
\begin{center}
\ \psfig{file=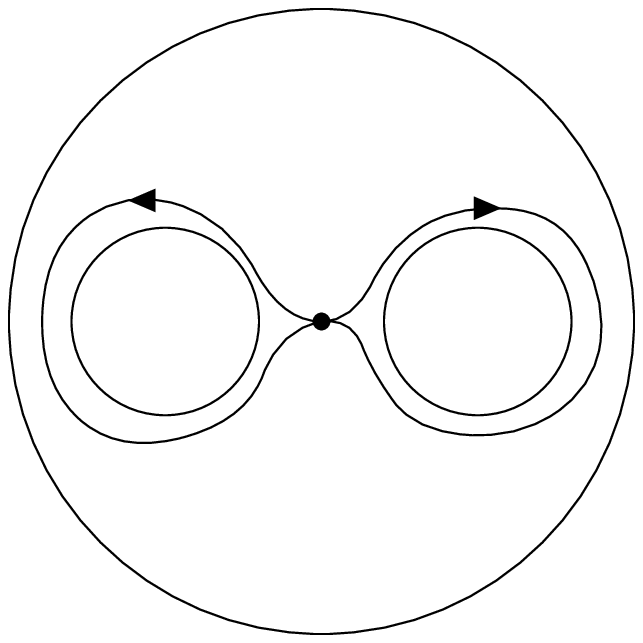,height=1.3truein}
\caption{Pants.}
\label{pants}
\end{center}
  \setlength{\unitlength}{1in}
  \begin{picture}(0,0)(0,0)
    \put(2.55,1.6){$a$}
    \put(3.35,1.6){$b$}
    \put(2.97,1.1){$x$}
  \end{picture}
\end{figure}

Choose a basepoint $x$ and pair of generators $a$ and $b$ for $\pi_{1}(P,x)$ as shown in the Figure \ref{pants} (so $\pi_{1}(P,x) \cong F(a,b)$, the free group on $a$ and $b$).
As in the closed case, we denote the set of essential closed curves in $P$ by ${\cal C}rv(P)$.
Any element of ${\cal C}rv(P)$ may then be represented by a cyclically reduced word, unique up to cyclic permutation and taking inverses (see e.g. \cite{LS}).
We write such an element as $a^{x_{1}}b^{y_{1}}...a^{x_{n}}b^{y_{n}}$ with $x_{j} \neq 0$ and $y_{j} \neq 0$ for each $j = 1,...,n$, except possibly when $n=1$, in which case we write the element as $a^{x_{1}}$, $b^{y_{1}}$, or $a^{x_{1}}b^{y_{1}}$.

The geometric intersection number of $\gamma \in {\cal C}rv(P)$ with $l_{1}$ is given by
$$i(\gamma,l_{1}) = \sum_{j=1}^{n}|x_{j}|.$$
where $\gamma$ is given as $a^{x_{1}}b^{y_{1}}...a^{x_{n}}b^{y_{n}}$ as above.
To see this, we construct a cover $\pi: \widetilde{P} \rightarrow P$ using cut and paste techniques so that
\begin{itemize}
\item $\gamma$ lifts to a curve $\widetilde{\gamma}$ in this cover
\item $\widetilde{\gamma}$ can be homotoped so that every component $\widetilde{l_{1}}$ of $\pi^{-1}(l_{1})$ intersects $\widetilde{\gamma}$ at most once.
\end{itemize}
It follows that
$$i(\gamma,l_{1}) = i(\widetilde{\gamma},\pi^{-1}(l_{1})) = \sum_{j=1}^{n}|x_{j}|.$$

Similarly, 
$$i(\gamma,l_{2}) = \sum_{j=1}^{n}|y_{j}|.$$

A similar type of argument can be used to show that
$$i(\gamma,w_{1}) = 2n$$
with the exception that
$$i(a^{x_{1}},w_{1}) = i(b^{y_{1}},w_{1}) = 0.$$
Formulae for $i(\gamma,l_{3})$ and $i(\gamma,w_{i})$, for $i=2,3$, can be obtained by choosing a different basis for $\pi_{1}(P,x)$.\\

Now suppose that the pair of pants, $P$, is an incompressible subsurface of a closed surface $S$.
Fixing a hyperbolic structure on $S$, we may assume that $P$ has geodesic boundary.
If $\alpha \in {\cal S}(S)$, and $\gamma,\gamma' \in {\cal C}(S)$, and $\gamma$ and $\gamma$ both have representatives lying entirely in $P$, then the geodesic representatives of $\gamma$ and $\gamma'$ also lie in $P$.
Taking the geodesic representative for $\alpha$, we see that $\alpha \cap P$ is a finite union of arcs of types $l_{1},l_{2},l_{3},w_{1},w_{2},w_{3}$.
It follows that if $i(\gamma,l_{i}) = i(\gamma',l_{i})$ and $i(\gamma,w_{i}) = i(\gamma',w_{i})$, in $P$, for each $i=1,2,3$, then $i(\gamma,\alpha) = i(\gamma',\alpha)$.
This proves the following.

\begin{proposition} \label{curvesinpants}
If $\gamma,\gamma' \in {\cal C}rv(S)$ have representatives contained in an incompressible pair of pants, $P \subset S$, then $\gamma \equiv_{si} \gamma'$ if and only if $\gamma$ and $\gamma'$ have the same geometric intersection number with each of the $6$ essential arcs in $P$.
\end{proposition}
\rightline{$\Box$}

\subsection{Counterexamples}

Before we construct our counterexamples, we will need the following.

\begin{proposition} \label{pushingforward}
Suppose $\gamma,\gamma' \in {\cal C}rv(S)$ and $\pi : \widetilde{S} \rightarrow S$ is a finite sheeted cover in which both $\gamma$ and $\gamma'$ lift to curves $\widetilde{\gamma}$ and $\widetilde{\gamma}'$ respectively.
If $\widetilde{\gamma} \equiv_{si} \widetilde{\gamma}'$, then $\gamma \equiv_{si} \gamma'$.
\end{proposition}

\noindent
{\bf Proof. }  We assume $\widetilde{\gamma} \equiv_{si} \widetilde{\gamma}'$ and let $\alpha \in {\cal S}(S)$ be arbitrary.
We write $\widetilde{\alpha} = \pi^{-1}(\alpha)$, and note that
$$\widetilde{\alpha} = \coprod_{j=1}^{k} \widetilde{\alpha}_{j}$$
where $\widetilde{\alpha}_{j} \in {\cal S}(\widetilde{S})$ for each $j=1,...,k$.
We then have
$$i(\alpha,\gamma) = i(\widetilde{\alpha},\widetilde{\gamma}) = \sum_{j=1}^{k}i(\widetilde{\alpha}_{j},\widetilde{\gamma}) = \sum_{j=1}^{k}(\widetilde{\alpha}_{j},\widetilde{\gamma}') = i(\widetilde{\alpha},\widetilde{\gamma}') = i(\alpha,\gamma')$$
\rightline{$\Box$}

\vspace{.3cm}

\begin{figure}[htb]
\begin{center}
\ \psfig{file=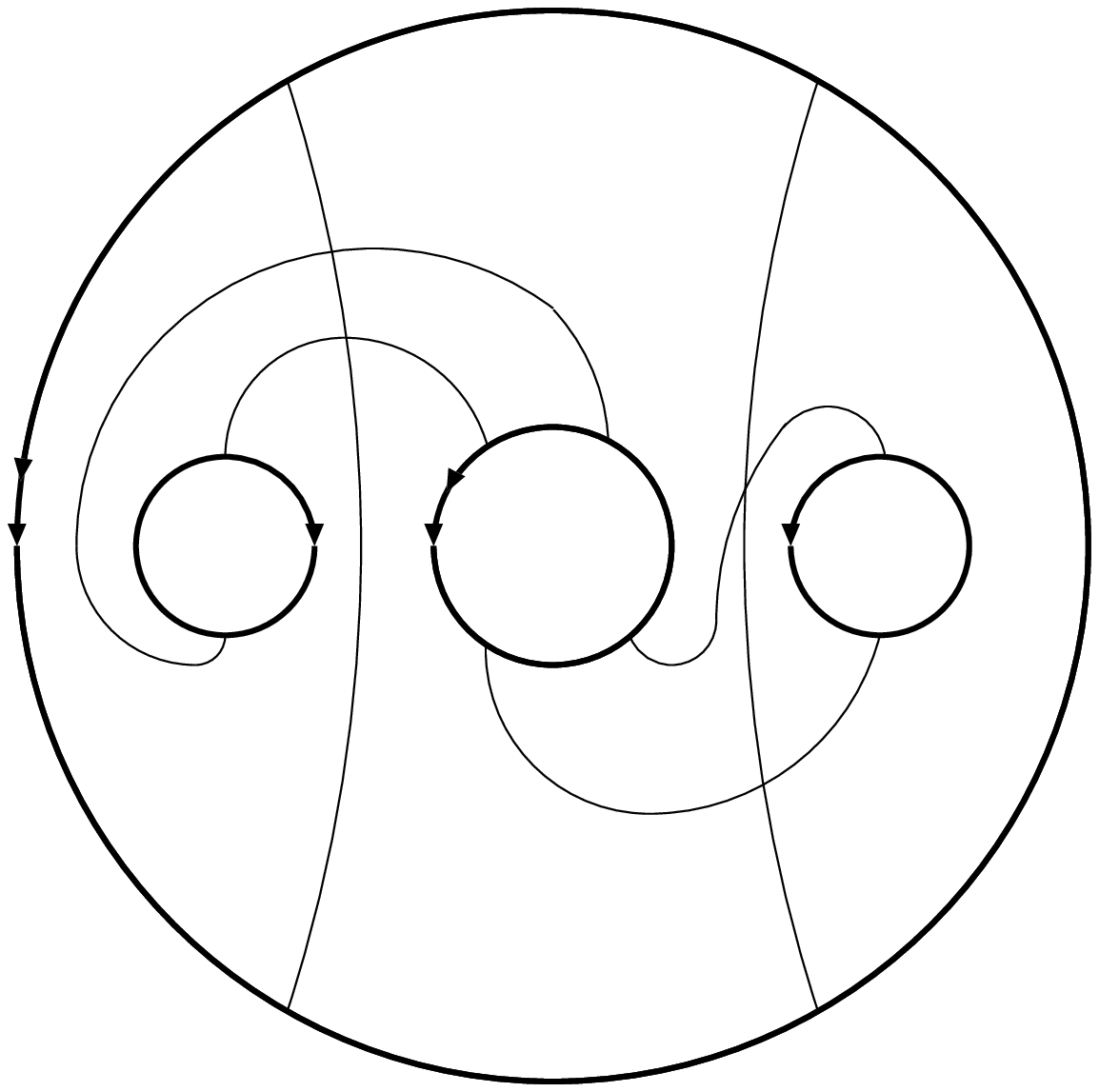,height=2truein}
\caption{The immersed curve $\gamma_{0}$.}
\label{immersedcurve}
\end{center}
  \setlength{\unitlength}{1in}
  \begin{picture}(0,0)(0,0)
    \put(2.4,1.6){$A$}
    \put(2.95,1.6){$B$}
    \put(3.53,1.6){$A'$}
    \put(3.8,2.3){$B'$}
    \put(2.9,2.2){$\gamma_{0}$}
  \end{picture}
\end{figure}

To construct our counterexample, we begin with a surface $S$ of genus $2$ and the curve $\gamma_{0} \in {\cal C}rv(S)$ shown in Figure \ref{immersedcurve}.
The surface $S$ is obtained by considering the sphere with $4$ holes shown and gluing boundary components $B$ to $B'$ by a dilation and $A$ to $A'$ by a reflection through the vertical line cutting the figure in half.
In the surface $S$, we refer to these curves as $B$ and $A$ respectively.
The arcs then match up to give the closed curve we call $\gamma_{0} \in {\cal C}rv(S)$.
One can check that cutting $S$ along $\gamma_{0}$ gives $2$ octagons, and we provide $S$ a hyperbolic metric so that each is a regular all-right octagon.
The full preimage of $\gamma_{0}$ in ${\mathbb H}^{2}$ (the universal cover) gives a tessellation by regular all-right octagons, so that in particular, every complete geodesic in ${\mathbb H}^{2}$ transversely intersects $\gamma_{0}$.
Therefore, every geodesic in $S$ transversely intersects $\gamma_{0}$, and $\gamma_{0}$ fills $S$.

\begin{figure}[htb]
\begin{center}
\ \psfig{file=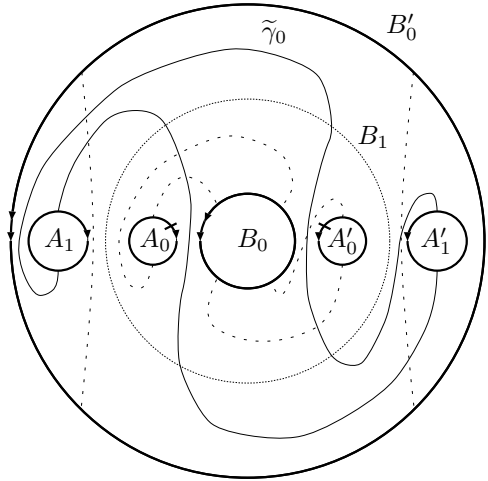,height=2.5truein}
\caption{A $2$-fold cover of $S$.}
\label{2foldcover}
\end{center}
  \setlength{\unitlength}{1in}
  \begin{picture}(0,0)(0,0)
    \put(1.96,1.83){$A_{1}$}
    \put(2.46,1.83){$A_{0}$}
    \put(2.97,1.83){$B_{0}$}
    \put(3.44,1.83){$A_{0}'$}
    \put(3.93,1.83){$A_{1}'$}
    \put(3.6,2.37){$B_{1}$}
    \put(3.75,2.95){$B_{0}'$}
    \put(3.1,2.92){$\widetilde{\gamma}_{0}$}
  \end{picture}
\end{figure}

Now consider the $2$-fold covering space $\pi:\widetilde{S} \rightarrow S$ corresponding to the kernel of the homomorphism induced by mod $2$ intersection number with $B$.
One easily checks that $\gamma_{0}$ has $0$ mod $2$ intersection number with $B$, and hence has two distinct lifts to $\widetilde{S}$.
Figure \ref{2foldcover} shows this two fold cover (the gluings are as indicated).
$A$ and $B$ both lift in this cover (one of the lifts of $B$, labeled $B_{1}$, is drawn as a dotted circle to help clarify the picture).
Both lifts of $\gamma_{0}$ are shown; one is drawn solid, which we call $\widetilde{\gamma}_{0}$, and the other is dashed.

\begin{figure}[htb]
\begin{center}
\ \psfig{file=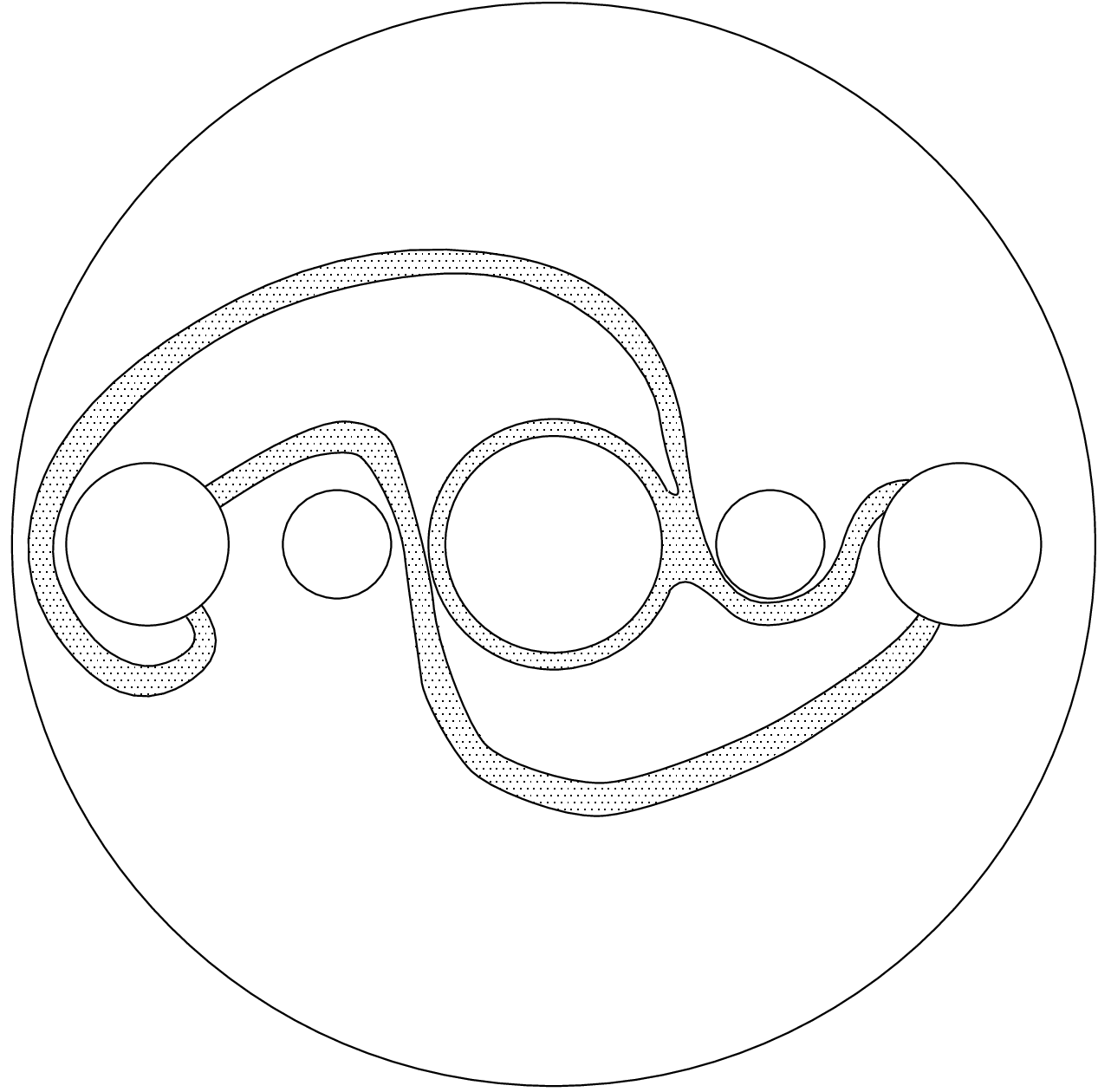,height=2.5truein}
\caption{Pants in $\widetilde{S}$.}
\label{pantsinS}
\end{center}
\end{figure}

Now we consider an incompressible pair of pants $P \subset \widetilde{S}$ with two components of $\partial P$ being $\widetilde{\gamma}_{0}$ and $B_{0}$.
We take $P$ as a regular neighborhood of $\widetilde{\gamma}_{0} \cup B_{0} \cup \rho$, where $\rho$ is a short arc connecting $\widetilde{\gamma}_{0}$ to $B_{0}$.
In Figure \ref{pantsinS} we have redrawn $\widetilde{S}$ indicating $P$ as the shaded subsurface.

We orient $\widetilde{\gamma}_{0}$ and $B_{0}$, pick a point $x \in int(P)$, and choose a basis $\{ a,b \}$ for $\pi_{1}(P,x)$ as in the previous section, so that $a$ and $b$ are freely homotopic (as oriented curves) to $\widetilde{\gamma}_{0}$ and $B_{0}$, respectively.

Now, let $\widetilde{\gamma}, \widetilde{\gamma}' \in {\cal C}rv(\widetilde{S})$ be curves contained in $P$ given by the following words in $a$ and $b$.
$$\widetilde{\gamma} = a^{-2}bab^{-1} \mbox{ , } \widetilde{\gamma}' = a^{-2}bab$$
and set
$$\gamma = \pi(\widetilde{\gamma}) \mbox{ , } \gamma' = \pi(\widetilde{\gamma}').$$

Note that $\widetilde{\gamma}$ and $\widetilde{\gamma}'$ are lifts of $\gamma$ and $\gamma'$.
For if this were not the case, then $\widetilde{\gamma}$, say, would be a $2:1$ cover of $\gamma$ in $\widetilde{S}$, and the covering transformation would have to take $\widetilde{\gamma}$ onto itself.
This would imply that the covering map restricts to a self map of $P$, which is impossible since there are no fixed point free orientation preserving involutions of a pair of pants.

\vspace{.3cm}

\noindent
{\bf Claim 1. }  Both $\gamma$ and $\gamma'$ fill $S$.

\noindent
{\bf Proof. }  This follows from the fact that $\gamma_{0}$ fills $S$.
To see this, consider our pair of pants $P \subset \widetilde{S}$ with its hyperbolic structure (with geodesic boundary), and represent all curves by geodesics.
Since neither of $\widetilde{\gamma}$ nor $\widetilde{\gamma}'$ is a multiple of a boundary component of $P$, it is easy to see that any curve which transversely intersects $\widetilde{\gamma}_{0}$ must transversely intersect both $\widetilde{\gamma}$ and $\widetilde{\gamma}'$.

Now any complete geodesic $\sigma$ in $S$ transversely intersects $\gamma_{0}$ since $\gamma_{0}$ fills $S$.
Any such point of intersection lifts to a point of intersection of $\widetilde{\gamma_{0}}$ with the preimage of $\sigma$.
This gives rise to points of transverse intersection of both $\widetilde{\gamma}$ and $\widetilde{\gamma}'$ with this lift of $\sigma$, by the remarks of the previous paragraph.
These points push down to points of intersection of $\sigma$ with $\gamma$ and $\gamma'$.
Therefore both $\gamma$ and $\gamma'$ fill $S$.\\

\noindent
{\bf Claim 2. }  $\gamma \equiv_{si} \gamma'$.

\noindent
{\bf Proof. }  We show that $\widetilde{\gamma} \equiv_{si} \widetilde{\gamma}'$ and apply Proposition \ref{pushingforward} to prove this claim.
By Proposition \ref{curvesinpants}, we need only check that $i(\widetilde{\gamma},l_{i}) = i(\widetilde{\gamma}',l_{i})$ and $i(\widetilde{\gamma},w_{i}) = i(\widetilde{\gamma}',w_{i})$, for $i=1,2,3$.
Using the method described in the previous subsection, one easily obtains
$$i(\widetilde{\gamma},l_{1}) = 3 = i(\widetilde{\gamma}',l_{1}) \mbox{ , } i(\widetilde{\gamma},l_{2}) = 2 = i(\widetilde{\gamma}',l_{2}) \mbox{ , } i(\widetilde{\gamma},l_{3}) = 3 =i(\widetilde{\gamma}',l_{3})$$
$$i(\widetilde{\gamma},w_{1}) = 4 = i(\widetilde{\gamma}',w_{1}) \mbox{ , } i(\widetilde{\gamma},w_{2}) = 4 = i(\widetilde{\gamma}',w_{2}) \mbox{ , } i(\widetilde{\gamma},w_{3}) = 6 =i(\widetilde{\gamma}',w_{3})$$

\vspace{.3cm}

\noindent
{\bf Claim 3. }  $\gamma \not\equiv_{h} \gamma'$.

\noindent
{\bf Proof. }  By Corollary \ref{h2hom}, it suffices to show that $[\gamma] \neq \pm [\gamma'] \in H_{1}(S;{\mathbb Z})$.
We have
$$[\gamma] = \pi_{*}([\widetilde{\gamma}]) = \pi_{*}(-2[\widetilde{\gamma}_{0}]+[B_{0}]+[\widetilde{\gamma}_{0}]-[B_{0}])
=\pi_{*}(-[\widetilde{\gamma}_{0}]) = -[\gamma_{0}]$$
$$[\gamma'] =\pi_{*}([\widetilde{\gamma}']) = \pi_{*}(-2[\widetilde{\gamma}_{0}]+[B_{0}]+[\widetilde{\gamma}_{0}]+[B_{0}])
=\pi_{*}(-[\widetilde{\gamma}_{0}]+2[B_{0}]) = -[\gamma_{0}]+2[B]$$ 
So, if $[\gamma] = \pm [\gamma']$ then we obtain
$$-[\gamma_{0}] = -[\gamma_{0}] + 2[B] \Leftrightarrow 2[B] = 0$$
or
$$-[\gamma_{0}] = [\gamma_{0}] - 2[B] \Leftrightarrow [B] = [\gamma_{0}].$$
Since $[B] \neq 0$, and since $[B]$ cannot possibly be $[\gamma_{0}]$ ($[\gamma_{0}]$ has non-trivial algebraic intersection number with $[A]$, for example) both of these give a contradiction, hence $[\gamma] \neq \pm [\gamma'] \in H_{1}(S;{\mathbb Z})$ and $\gamma \not\equiv_{h} \gamma'$.
In a similar fashion, one can construct many other examples.

We have proved

\begin{lemma} \label{sinot2h}
There exists curves $\gamma,\gamma' \in {\cal C}rv(S)$, each of which fill $S$, such that $\gamma \equiv_{si} \gamma'$, yet $\gamma \not\equiv_{h} \gamma'$.
\end{lemma}
\rightline{$\Box$}

Proposition \ref{tr2h2tr}, Proposition \ref{bf2si2bf}, Corollary \ref{h2si}, and Lemma \ref{sinot2h} combine to prove Theorem \ref{main}.

\section{Alternatives and related questions}  \label{theend}

Suppose $\gamma, \gamma' \in {\cal C}rv(S)$ and $\gamma \equiv_{h} \gamma'$.
Let $f:S \rightarrow T^{2}$ be a continuous function from $S$ to the torus $T^{2}$.
It follows from Corollary \ref{h2hom} that the map $f_{*}:\pi_{1}(S) \rightarrow \pi_{1}(T^{2}) \cong {\mathbb Z}^{2}$ must have $f_{*}(\gamma) = f_{*}(\gamma')$.
On a torus this implies $f(\gamma) = f(\gamma')$ in ${\cal C}rv(T^{2})$, and in particular, $f(\gamma) \equiv_{si} f(\gamma')$ ($\equiv_{si}$ is defined as in the hyperbolic case).

Now, let $S'$ be any orientable surface of genus $g > 1$, and let $f:S \rightarrow S'$ be any continuous function.
By Corollary \ref{h2maphyper}, $f(\gamma) \equiv_{h} f(\gamma')$ and so by Theorem \ref{main}, $f(\gamma) \equiv_{si} f(\gamma')$.
This leads us to make the following

\begin{definition}
Given $\gamma, \gamma' \in {\cal C}rv(S)$, say that $\gamma$ and $\gamma'$ are {\em strongly simple intersection equivalent}, $\gamma \equiv_{ssi} \gamma'$, if for every orientable surface $S'$ and continuous map $f:S \rightarrow S'$, we have $f(\gamma) \equiv_{si} f(\gamma')$.
\end{definition}

The preceding arguments prove

\begin{theorem}
If $\gamma, \gamma' \in {\cal C}rv(S)$ and $\gamma \equiv_{h} \gamma'$, then $\gamma \equiv_{ssi} \gamma'$.
\end{theorem}
\rightline{$\Box$}

Strong simple intersection equivalence is certainly stronger than simple intersection equivalence as the counterexamples in the previous section show.

\begin{question}
Is is true that if $\gamma, \gamma' \in {\cal C}rv(S)$ have $\gamma \equiv_{ssi} \gamma'$, then $\gamma \equiv_{h} \gamma'$?
\end{question}

It is not clear what the answer to this question should be.
The initial problem encountered is that there seems to be no real way to check that a pair of curves are strongly simple intersection equivalent without assuming that they are hyperbolically equivalent.

In particular, we pose the following

\begin{question}
Given a pair of curves, $\gamma, \gamma' \in {\cal C}rv(S)$, is there an algorithm to decide whether or not $\gamma \equiv_{ssi} \gamma'$?
\end{question}

We have not touched on any other possible characterizations of hyperbolic equivalence, and we have also not mentioned anything about the higher dimensional case (in particular, the $3$-dimensional case \cite{Mas}).
For more on this see \cite{Sch}, \cite{GR}, and particularly Conjecture 4.1 of \cite{And}.

\noindent
Address:\\
Department of Mathematics\\
Barnard College at Columbia University\\
2990 Broadway MC 4448\\
New York, NY 10027-6902\\
Phone: (512)-854-9235\\
email: clein@math.columbia.edu\\

\end{document}